\documentclass[a4paper,11pt,fleqn]{article}
\usepackage{amsmath,amssymb,amsthm,graphicx,subfigure,float,caption,epstopdf,tabularx,color, bm,amsfonts,epic}
\usepackage[top=1in, bottom=1in, left=1.25in, right=1.25in]{geometry}
\usepackage{appendix}
\usepackage{diagbox}
\allowdisplaybreaks
\newtheorem{thm}{Theorem}[section]
\newtheorem{lem}{Lemma}[section]
\newtheorem{coro}{Corollary}[section]
\newtheorem{den}{Definition}[section]

\newtheorem{rmk}{Remark}[section]
\newtheorem*{prf}{Proof}
\numberwithin{equation}{section}

\usepackage{indentfirst}
\graphicspath{{Fig}}

\begin{document}
\title{A linearly implicit structure-preserving Fourier pseudo-spectral scheme for the damped nonlinear Schr\"{o}dinger
equation in three dimensions}
\author{Chaolong Jiang$^{2}$, \quad Yongzhong Song$^{1}$, \quad Yushun Wang$^{1,}$\footnote{Correspondence author. Email:
wangyushun@njnu.edu.cn.}\\
{\small $^1$ Jiangsu Key Laboratory for Numerical Simulation of Large Scale Complex Systems,}\\
{\small School of Mathematical Sciences,  Nanjing Normal University,}\\
{\small  Nanjing 210023, China}\\
{\small $^2$ School of Statistics and Mathematics, }\\
{\small Yunnan University of Finance and Economics, Kunming 650221, China}\\
}
\date{}
\maketitle

\begin{abstract}
In this paper, we propose a linearly implicit Fourier pseudo-spectral scheme,
which preserves the total mass and energy conservation laws for
the damped nonlinear Schr\"{o}dinger equation in three dimensions.
With the aid of the semi-norm equivalence between
the Fourier pseudo-spectral method and the finite difference
method, an optimal $L^{2}$-error estimate for the proposed method without any restriction on the grid ratio is established by analyzing the real and imaginary parts of the error function. Numerical results are addressed to confirm our theoretical analysis. \\[2ex]
\end{abstract}

\section{Introduction}
The damped nonlinear Schr\"{o}dinger (DNLS) equation arises in various fields of physics, such as water
waves \cite{SHC05}, nonlinear optics \cite{HK95}, and plasma physics \cite{GRH80}.
In this paper, we consider the following DNLS equation in three dimensions (3D)
\begin{align}\label{3DNLD_eq:1.1}
 \text{i}\psi_{t}(t,{\bm x})+\Delta \psi(t,{\bm x})+\beta \big|\psi(t,{\bm x})\big|^{2}\psi(t,{\bm x})+\text{i}\gamma \psi(t,{\bm x})=0,\ {\bm x}\in\Omega,\ 0<t\leq T,
\end{align}
with $(l_{1},l_{2},l_{3})$-periodic boundary conditions
\begin{align*}
&\psi(t,x,y,z)=\psi(t,x+l_{1},y,z),\ \psi(t,x,y,z)=\psi(t,x,y+l_{2},z),\nonumber\\
& \psi(t,x,y,z)=\psi(t,x,y,z+l_{3}),\ (x,y,z)\in\Omega,\ 0<t\leq T,
\end{align*}
and initial condition
\begin{align*}
\psi(0,{\bm x})=\psi_{0}({\bm x}),\ {\bm x}\in\Omega,
\end{align*}
where $\text{i}=\sqrt{-1}$ is the complex unit, $t$ is the time variable, ${\bm x}=(x,y,z)$ is the spatial variable, $\psi:=\psi(t,{\bm x})$ is the complex-valued wave function, $\Delta =\partial_{xx}+\partial_{yy}+\partial_{zz}$ is the usual Laplace operator, $\beta$ is a
given real constant, $\gamma> 0$ represents dissipation, $\Omega=[0,l_{1}]\times[0,l_{2}]\times [0,l_{3}]\subset\mathbb{R}^{3}$
 and $\psi_{0}:=\psi_{0}({\bm x})$ is a given $(l_{1},l_{2},l_{3})$-periodic complex-valued function.
When $\gamma=0$, the DNLS equation \eqref{3DNLD_eq:1.1} reduces to the classical NLS equation.

Let $u:=u(t,{\bm x})=e^{\gamma t}\psi(t,{\bm x})$. Eq. \eqref{3DNLD_eq:1.1} can be rewritten as
\begin{align}\label{3DNLD-eq-1.2}
 &\text{i}u_{t}(t,{\bm x})+ \Delta u(t,{\bm x})+\beta e^{-2\gamma t}\big|u(t,{\bm x})\big|^{2}u(t,{\bm x})=0,\ {\bm x}\in\Omega,\ 0<t\leq T,
\end{align}
with the periodic conditions
\begin{align*}
&u(t,x,y,z)=u(t,x+l_{1},y,z),\ u(t,x,y,z)=u(t,x,y+l_{2},z),\\
 & u(t,x,y,z)=u(t,x,y,z+l_{3}),\ (x,y,z)\in\Omega,\ 0<t\leq T,
\end{align*}
and initial condition
\begin{align*}
u(0,{\bm x})=\psi_{0}({\bm x}),\ {\bm x}\in\Omega.
\end{align*}
The DNLS \eqref{3DNLD-eq-1.2} admits the following mass conservation law
\begin{align*}
M(t):=\int_{\Omega}|u(t,{\bm x})|^2d{\bm x}\equiv \int_{\Omega}|u(0,{\bm x})|^2d{\bm x}=\int_{\Omega}|\psi_0({\bm x})|^2d{\bm x}:=M(0),\ 0<t\leq T,
\end{align*}
 and energy conservation law
\begin{align*}
E(t):=&\int_{\Omega}|\nabla u(t,{\bm x})|^2d{\bm x}-\frac{\beta}{2}e^{-2\gamma t}\int_{\Omega}|u(t,{\bm x})|^4d{\bm x}\\
&~~~~~~~~~~~~~~~~~~~~~~~~~-\gamma\beta\int_{0}^{t}e^{-2\gamma \nu}\int_{\Omega}|u(\nu,{\bm x})|^4d{\bm x}d\nu\equiv E(0),\ 0<t\leq T.
\end{align*}

The analysis and numerical solution of the DNLS equation has been widely investigated.
Fibich \cite{GF01} analyzed the effect of linear damping (absorption) on critical self-focusing NLS equation.
Tsutsumi \cite{MT84} studied the global solutions of the DNLS equation.
The regularity and existence of attractors for a weakly DNLS equation were investigated in Ref. \cite{Temam12}.
Known strategies to solve the DNLS equation numerically include the time-splitting sine pseudo-spectral method \cite{ABB13,BC13b,BJ03}, finite difference 
methods \cite{BM2018,DFP81,FZQ16,HDY17,IJB81,MNS13,Peranich81}, and the finite element method \cite{ADKM98}. However, most of existing schemes are fully implicit and only considered for the DNLS equation in one or two dimensions. For the fully implicit schemes, one needs to solve a system of nonlinear equations at each time step, which leads to expensive costs. Therefor, how to design highly efficient numerical schemes for the DNLS equation in 3D attracts a lot of interest.

Error estimates for different numerical methods of the DNLS equation in one dimension have been established.
For the Fourier spectral Galerkin method, we refer to Ref. \cite{xiang98} for details. For the finite difference method, we refer to Ref. \cite{ZL01} for details.
Unconditionally optimal error analyses of conservative methods for the
DNLS equation in one dimension were conducted in Refs. \cite{JCW17,ZWZ11}.
In fact, their proofs for conservation schemes \cite{JCW17,ZWZ11} rely heavily on not only the discrete conservative property
 but also the discrete version of the Sobolev inequality in one dimension
\begin{align*}
\|f\|_{L^{\infty}}\leq C\|f\|_{H^{1}},\ \forall f\in H^{1}(\Omega)\  \text{with\ the bounded\ domain}\ \Omega\subset \mathbb{R},
\end{align*}
which immediately implies an a priori uniform bound for $||f||_{L^\infty}$.
However, the extension of the discrete version of the above Sobolev inequality is no longer valid in 3D.
Thus the techniques used in Refs. \cite{JCW17,ZWZ11} cannot be extended directly to 3D.
Due to the difficulty in obtaining such a priori bound for the numerical solution,
few error estimates are obtained in the literature for the DNLS equation in 3D. With the aid of the
classical inverse inequality, Zhang \cite{zhang04} established an optimal $L^2$-error estimate for a fully implicit finite difference (IFD) scheme of the DNLS equation in 3D. However, the result requires a strict restriction on the grid ratio. In this paper, we fucus on establishing an optimal error estimate without any restrictions on the grid ratio for numerical schemes of the DNLS equation in 3D.

It is well-known that the Fourier spectral method has been a powerful tool to solve partial differential equations (PDEs), both theoretically and numerically \cite{Boyd01,ST06}. With the advent of the fast Fourier transform (FFT),
the Fourier spectral method provides a numerical discretization with the convergency of so-called infinite order and high efficiency. Actually, for high dimensional problems with periodic boundary conditions, the efficiency of the Fourier spectral method is comparable to that of the finite difference method. Bridges and Reich \cite{BR01b} first introduced the idea of Fourier spectral discretization to
construct a multi-symplectic integrator for Hamiltonian PDEs. Motivated
by the theory of Bridges and Reich, Chen and Qin \cite{CQ01} proposed a symplectic and
multi-symplectic Fourier pseudo-spectral method for the nonlinear Schr\"odinger equation with periodic boundary conditions. Later, different kinds of structure-preserving Fourier pseudo-spectral schemes were developed (e.g., see Refs. \cite{CSZ11,GCW14,KZCDH10,LW15}). However,
to our best knowledge, most existing works on structure-preserving Fourier pseudo-spectral schemes up to now focus on Hamiltonian PDEs and fail to damped cases.  On the other hand, apart from numerical implementations, theoretical analysis on the convergence of the structure-preserving Fourier pseudo-spectral schemes are still highly desired. In Ref. \cite{GWWC17}, Gong et al. established an optimal $L^2$-error estimate for
a conservative Fourier pseudo-spectral scheme for the nonlinear Schr\"odinger equation in 2D.
 The scheme and analysis technique can be generalized to the damped nonlinear Schr\"odinger equation in 2D, but not for the case in 3D.
 This is due to the fact that the discrete version of interpolation
 inequalities \cite{BC13,WGX13} in 2D is not valid for the 3D case.
In Refs. \cite{CHWG15,CWG16,JCWL17}, the authors established an optimal error estimate, without any restriction on the grid ratio, for an
energy-preserving method of the 3D Maxwell's equations in discrete $L^2$-norm, which is however only
useful for linear problems. Thus, the standard convergence analysis of structure-preserving Fourier pseudo-spectral methods for damped Hamiltonian PDEs in 3D is still on the early stage.

In this paper, we propose a linearly implicit and conservative Fourier pseudo-spectral (LI-CFP) scheme for the DNLS equation in 3D by using in time a linearly implicit energy-preserving method \cite{ABB13} and in space the standard Fourier pseudo-spectral method. For the linearly implicit scheme, we only require to solve a linear system of equations at each time step, which leads to considerably lower costs than the implicit ones. Using a new introduced semi-norm equivalence between the Fourier pseudo-spectral method and the finite difference method, the projection and interpolation theories, together with a linearized technique for the finite difference method \cite{SW17}, we show that the proposed scheme is unconditionally convergent with the order of $O(N^{-s}+\tau^2)$ in discrete $L^2$-norm, where $N$ is the number of collocation points used in the spectral method and $\tau$ is the time step.

The outline of this paper is
organized as follows. In Section 2, a linearly implicit Fourier pseudo-spectral scheme, which can preserve the discrete mass and energy conservation laws, for the DNLS equation in 3D and some important Lemmas are presented. In Section 3, we prove that the numerical scheme is uniquely solvable.
In Section 4, an a priori estimate for the proposed
scheme is established in discrete $L^{2}$-norm.
Some numerical experiments are presented in Section 5. We draw some conclusions in Section 6.

\section{Construction of the linearly implicit and conservative scheme}
Let $\Omega_{h}=\{(x_{j_1},y_{j_2},z_{j_3})|x_{j_{1}}=j_{1}h_{1},y_{j_{2}}=j_{2}h_{2},z_{j_{3}}=j_{3}h_{3};\ 0\leq j_{r}\leq N_{r}-1, r=1,2,3\}$ be a partition of $\Omega$
with the grid size $h_{r}=\frac{l_{r}}{N_{r}}$, where $N_{r}$ is an even number. Denote $h=\text{max}\{h_1,h_2,h_3\}$. Let $\Omega_{\tau}=\{t_{n}|t_{n}=n\tau; 0\leq n\leq M\}$
be a uniform partition of $[0,T]$ with the time step $\tau=T/M$, $\Omega_{h\tau}=\Omega_{h}\times\Omega_{\tau}$, and denote
\begin{align*}
J_{h}=\{\vec{j}=(j_1,j_2,j_3)|0\leq j_{r}\leq N_{r}-1, r=1,2,3\}.
\end{align*}
A discrete mesh function ${ u}_{\vec{j}},\ \vec{j}\in J_{h}$ defined on $\Omega_h$ is said to satisfy the periodic boundary conditions if and only if
\begin{align}\label{DNLS-PBS}
&\left\lbrace
  \begin{aligned}
&x-\text{periodic}:\ u_{j_1,j_2,j_3}=u_{j_1+N_1,j_2,j_3},\ j_r=0,1,2\cdots,N_r-1,\ r=2,3,\\
&y-\text{periodic}:\ u_{j_1,j_2,j_3}=u_{j_1,j_2+N_2,j_3},\ j_r=0,1,2\cdots,N_r-1,\ r=1,3,\\
&z-\text{periodic}:\ u_{j_1,j_2,j_3}=u_{j_1,j_2,j_3+N_3},\ j_r=0,1,2\cdots,N_r-1,\ r=1,2.
\end{aligned}\right.\
\end{align}
Denoting $\{U_{\vec{j}}^{n}|\vec{j}\in J_h\}$ and $\{V_{\vec{j}}^{n}|\vec{j}\in J_h\}$ as two grid functions defined on $\Omega_{h\tau}$.
We introduce the following notations:
\begin{align*}
&\delta_{x}^{+} U_{\vec{j}}^{n}=\frac{U_{j_{1}+1,j_{2},j_{3}}^{n}-U_{j_{1},j_{2},j_{3}}^{n}}{h_{1}},\
\delta_{y}^{+} U_{\vec{j}}^{n}=\frac{U_{j_{1},j_{2}+1,j_{3}}^{n}-U_{j_{1},j_{2},j_{3}}^{n}}{h_{2}},\\
&\delta_{z}^{+} U_{\vec{j}}^{n}=\frac{U_{j_{1},j_{2},j_{3}+1}^{n}-U_{j_{1},j_{2},j_{3}}^{n}}{h_{3}},\ \delta_{t} U_{\vec{j}}^{n}
=\frac{U_{\vec{j}}^{n+1}-U_{\vec{j}}^{n-1}}{2\tau},\\
& \delta_{t}^+ U_{\vec{j}}^{n}
=\frac{U_{\vec{j}}^{n+1}-U_{\vec{j}}^{n}}{\tau},\ \widehat{U}_{\vec{j}}^{n}=\frac{U_{\vec{j}}^{n+1}+U_{\vec{j}}^{n-1}}{2},\ {U}_{\vec{j}}^{n+\frac{1}{2}}=\frac{U_{\vec{j}}^{n+1}+U_{\vec{j}}^{n}}{2}.
\end{align*}

{Let
\begin{align*}
\mathbb{ V}_{h}:&=\big\{{\bm U}|{\bm U}=({\bm U}_{0,0},{\bm U}_{1,0},\cdots,{\bm U}_{N_{2}-1,0},{\bm U}_{0,1},{\bm U}_{1,1},\cdots,\nonumber\\
&~~~~~~~~~~~~~~~~~{\bm U}_{N_{2}-1,1},\cdots,{\bm U}_{0,N_3-1},{\bm U}_{1,N_3-1},\cdots,{\bm U}_{N_{2}-1,N_3-1})^T\big\}\nonumber\\
&=\big\{{\bm U}|{\bm U}=({\bm U}_{0},{\bm U}_{1},\cdots,{\bm U}_{N_{3}-1})^T\big\}
\end{align*}
be the space of grid functions defined on $\Omega_h$ that satisfy the periodic boundary condition \eqref{DNLS-PBS}, where
\begin{align*}
&{ {\bm U}}_{j_2,j_3}=(U_{0,j_2,j_3},\cdots,U_{N_{1}-1,j_2,j_3}),\\
 &{{\bm U}}_{j_3}=({\bm U}_{0,j_3},\cdots,{\bm U}_{N_{2}-1,j_3}),\ 0\leq j_r\leq N_r-1,r=2,3.
\end{align*}}
For any two grid functions ${\bm U},{\bm V}\in\mathbb{ V}_{h}$,
we define the discrete inner product and notions as, respectively,
\begin{align*}
&\langle {\bm U},{\bm V}\rangle_{h}=h_{\Delta}\sum_{\vec{j}\in J_{h}}U_{\vec{j}}\bar{V}_{\vec{j}},\ {\|{\bm U}\|_{h}^2=h_{\Delta}\sum_{\vec{j}\in J_{h}}|U_{\vec{j}}|^2},\\
& {|{\bm U}|_{1,h_1}^2=h_{\Delta}\sum_{\vec{j}\in J_{h}}|\delta_x^+U_{\vec{j}}|^2,\ |{\bm U}|_{1,h_2}^2=h_{\Delta}\sum_{\vec{j}\in J_{h}}|\delta_y^+U_{\vec{j}}|^2,} \\
&{|{\bm U}|_{1,h_3}^2=h_{\Delta}\sum_{\vec{j}\in J_{h}}|\delta_z^+U_{\vec{j}}|^2, |{\bm U}|_{1,h}^2=|{\bm U}|_{1,h_1}^2+|{\bm U}|_{1,h_2}^2+|{\bm U}|_{1,h_3}^2},\\
&||{\bm U}||_{h,p}^p=h_{\Delta}\sum_{\vec{j}\in J_{h}}|U_{\vec{j}}|^{p},\ \|{\bm U}\|_{h,\infty}=\max\limits_{\vec{j}\in J_{h}}|U_{\vec{j}}|,
\end{align*}
where $h_{\Delta}=h_1h_2h_3$, $\bar{V}_{\vec{j}}$ denotes the conjugate of $V_{\vec{j}}$, and {$|\cdot|$ is the absolute value of $\cdot$.} We note that $\|{\bm U}\|_{h}$ and $\|{\bm U}\|_{h,\infty}$ are norms and called $L^2$- and $L^{\infty}$-norms, respectively. In addition, we denote $`\cdot$' as the componentwise product of the vectors, that is,
\begin{align*}
{\bm U}\cdot {\bm V}=&\big(U_{0,0,0}V_{0,0,0},\cdots,U_{N_1-1,0,0}V_{N_1-1,0,0},\cdots,U_{0,N_2-1,N_3-1}V_{0,N_2-1,N_3-1},\\
&\cdots,U_{N_1-1,N_2-1,N_3-1}V_{N_1-1,N_2-1,N_3-1}\big)^{T}.
\end{align*}
For brevity, we denote ${\bm U}\cdot {\bm U}$ as ${\bm U}^2$.
\begin{den}\label{DNLS:den2.1} {In this paper, for any matrices ${\bm A}=(a_{j,k})_{p,q}$ and ${\bm B}=(b_{j,k})_{r,m}$, where $p,\ q,\ r,\ m$ are nonnegative integers, the Kronecker product ${\bm A}\otimes{\bm B}$ is a $pr\times qm$ block matrix defined by $$
{\bm A}\otimes{\bm B}=\left(\begin{array}{lllll}
a_{1,1}{\bm B}& a_{1,2}{\bm B}  &\cdots &a_{1,q}{\bm B} \\
              a_{2,1}{\bm B}& a_{2,2}{\bm B}  &\cdots &a_{2,q}{\bm B} \\
              \vdots &\vdots  &\ddots &\vdots\\
              a_{p,1}{\bm B}&  a_{p,2}{\bm B}& \cdots & a_{p,q}{\bm B}\\
\end{array}
\right).
$$ }
\end{den}

\begin{coro}\label{DNLS:cor:2.1} {According to the definition \ref{DNLS:den2.1}, we can show that, for any matrices ${\bm A}=(a_{j,k})_{p,m}$, ${\bm B}=(b_{j,k})_{m,l}$, ${\bm C}=(c_{j,k})_{r,q}$, ${\bm D}=(d_{j,k})_{q,s}$ and ${\bm E}=(c_{j,k})_{m_1,m_2}$, ${\bm F}=(d_{j,k})_{m_2,m_3}$, where $p,\ m,\ l,\ r,\ q,\ s,\ m_1,\ m_2,\ m_3$ are nonnegative integers, the Kronecker product $\otimes$ satisfies
\begin{align*}
({\bm A}\otimes {\bm C}\otimes{\bm E})({\bm B}\otimes {\bm D}\otimes{\bm F})={\bm A}{\bm B}\otimes{\bm C}{\bm D}\otimes{\bm E}{\bm F},\ ({\bm A}\otimes {\bm C}\otimes{\bm E})^{H}={\bm A}^H\otimes {\bm C}^H\otimes{\bm E}^H,
\end{align*}
where ${\bm W}^{H}$ represents the conjugate transpose matrix of ${\bm W}.$ }
\end{coro}
\begin{lem}\label{DNLS:lem:2.1} {For any matrices ${\bm A}=(a_{j,k})_{N_1,N_1}$, ${\bm B}=(b_{j,k})_{N_2,N_2}$, and ${\bm C}=(c_{j,k})_{N_3,N_3}$,
we have
\begin{align*}
({\bm C}\otimes{\bm B}\otimes{\bm A}){\bm U}=\left(\begin{array}{c}
\sum_{j_3=0}^{N_3-1}\sum_{j_2=0}^{N_2-1}\sum_{j_1=0}^{N_1-1}c_{0,j_3}b_{0,j_2}a_{0,j_1}{U}_{j_1,j_2,j_3}\\
\vdots\\
\sum_{j_3=0}^{N_3-1}\sum_{j_2=0}^{N_2-1}\sum_{j_1=0}^{N_1-1}c_{0,j_3}b_{0,j_2}a_{N_1-1,j_1}{U}_{j_1,j_2,j_3}\\
\sum_{j_3=0}^{N_3-1}\sum_{j_2=0}^{N_2-1}\sum_{j_1=0}^{N_1-1}c_{0,j_3}b_{1,j_2}a_{0,j_1}{U}_{j_1,j_2,j_3}\\
\vdots\\
\sum_{j_3=0}^{N_3-1}\sum_{j_2=0}^{N_2-1}\sum_{j_1=0}^{N_1-1}c_{0,j_3}b_{1,j_2}a_{N_1-1,j_1}{U}_{j_1,j_2,j_3}\\
\vdots\\
\sum_{j_3=0}^{N_3-1}\sum_{j_2=0}^{N_2-1}\sum_{j_1=0}^{N_1-1}c_{1,j_3}b_{0,j_2}a_{0,j_1}{U}_{j_1,j_2,j_3}\\
\vdots\\
\sum_{j_3=0}^{N_3-1}\sum_{j_2=0}^{N_2-1}\sum_{j_1=0}^{N_1-1}c_{1,j_3}b_{0,j_2}a_{N_1-1,j_1}{U}_{j_1,j_2,j_3}\\
\vdots\\
\sum_{j_3=0}^{N_3-1}\sum_{j_2=0}^{N_2-1}\sum_{j_1=0}^{N_1-1}c_{N_3-1,j_3}b_{N_2-1,j_2}a_{0,j_1}{U}_{j_1,j_2,j_3}\\
\vdots\\
\sum_{j_3=0}^{N_3-1}\sum_{j_2=0}^{N_2-1}\sum_{j_1=0}^{N_1-1}c_{N_3-1,j_3}b_{N_2-1,j_2}a_{N_1-1,j_1}{U}_{j_1,j_2,j_3}\\
 \end{array}
 \right),
\end{align*}
where ${\bm U}\in\mathbb{V}_h$. }
\end{lem}
\begin{prf} {According to the definition \ref{DNLS:den2.1}, we have
\begin{align}\label{DNLS:eq-2.2}
({\bm B}\otimes {\bm A}){\bm U}_{j_3}&=\left(\begin{array}{cccc}
b_{0,0}{\bm A}&\cdots &b_{0,N_2-1}{\bm A}\\
              b_{1,0}{\bm A}&\cdots &b_{1,N_2-1}{\bm A} \\
              \vdots  &\ddots &\vdots\\
              b_{N_2-1,0}{\bm A}& \cdots & b_{N_2-1,N_2-1}{\bm A}\\
\end{array}
\right)\left(\begin{array}{c}
{\bm U}_{0,j_3}\\
{\bm U}_{1,j_3}\\
\vdots\\
{\bm U}_{N_2-1,j_3}\\
 \end{array}
 \right)\nonumber\\
 &=\left(\begin{array}{c}
\sum_{j_2=0}^{N_2-1}b_{0,j_2}{\bm A}{\bm U}_{j_2,j_3}\\
\sum_{j_2=0}^{N_2-1}b_{1,j_2}{\bm A}{\bm U}_{j_2,j_3}\\
\vdots\\
\sum_{j_2=0}^{N_2-1}b_{N_2-1,j_2}{\bm A}{\bm U}_{j_2,j_3}\\
 \end{array}
 \right)\nonumber\\
 &=\left(\begin{array}{c}
\sum_{j_2=0}^{N_2-1}\sum_{j_1=0}^{N_1-1}b_{0,j_2}a_{0,j_1}{U}_{j_1,j_2,j_3}\\
\vdots\\
\sum_{j_2=0}^{N_2-1}\sum_{j_1=0}^{N_1-1}b_{0,j_2}a_{N_1-1,j_1}{U}_{j_1,j_2,j_3}\\
\sum_{j_2=0}^{N_2-1}\sum_{j_1=0}^{N_1-1}b_{1,j_2}a_{0,j_1}{U}_{j_1,j_2,j_3}\\
\vdots\\
\sum_{j_2=0}^{N_2-1}\sum_{j_1=0}^{N_1-1}b_{1,j_2}a_{N_1-1,j_1}{U}_{j_1,j_2,j_3}\\
\vdots\\
\sum_{j_2=0}^{N_2-1}\sum_{j_1=0}^{N_1-1}b_{N_2-1,j_2}a_{0,j_1}{U}_{j_1,j_2,j_3}\\
\vdots\\
\sum_{j_2=0}^{N_2-1}\sum_{j_1=0}^{N_1-1}b_{N_2-1,j_2}a_{N_1-1,j_1}{U}_{j_1,j_2,j_3}\\
 \end{array}
 \right),
\end{align}
where $ 0\le j_3\le N_3-1.$
Thus, it follows from \eqref{DNLS:eq-2.2} that
\begin{align*}
({\bm C}\otimes{\bm B}\otimes{\bm A}){\bm U}&=\left(\begin{array}{cccc}
c_{0,0}{\bm B}\otimes {\bm A}&\cdots &c_{0,N_3-1}{\bm B}\otimes {\bm A}\\
              c_{1,0}{\bm B}\otimes {\bm A}&\cdots &c_{1,N_3-1}{\bm B}\otimes {\bm A} \\
              \vdots  &\ddots &\vdots\\
              c_{N_3-1,0}{\bm B}\otimes {\bm A}& \cdots & c_{N_3-1,N_3-1}{\bm B}\otimes {\bm A}\\
\end{array}
\right)\left(\begin{array}{c}
{\bm U}_{0}\\
{\bm U}_{1}\\
\vdots\\
{\bm U}_{N_3-1}\\
 \end{array}
 \right)\nonumber\\
 &=\left(\begin{array}{c}
\sum_{j_3=0}^{N_3-1}c_{0,j_3}{\bm B}\otimes {\bm A}{\bm U}_{j_3}\\
\sum_{j_3=0}^{N_3-1}c_{1,j_3}{\bm B}\otimes {\bm A}{\bm U}_{j_3}\\
\vdots\\
\sum_{j_3=0}^{N_3-1}c_{N_3-1,j_3}{\bm B}\otimes {\bm A}{\bm U}_{j_3}\\
 \end{array}
 \right)\nonumber\\
 &=\left(\begin{array}{c}
\sum_{j_3=0}^{N_3-1}\sum_{j_2=0}^{N_2-1}\sum_{j_1=0}^{N_1-1}c_{0,j_3}b_{0,j_2}a_{0,j_1}{U}_{j_1,j_2,j_3}\\
\vdots\\
\sum_{j_3=0}^{N_3-1}\sum_{j_2=0}^{N_2-1}\sum_{j_1=0}^{N_1-1}c_{0,j_3}b_{0,j_2}a_{N_1-1,j_1}{U}_{j_1,j_2,j_3}\\
\sum_{j_3=0}^{N_3-1}\sum_{j_2=0}^{N_2-1}\sum_{j_1=0}^{N_1-1}c_{0,j_3}b_{1,j_2}a_{0,j_1}{U}_{j_1,j_2,j_3}\\
\vdots\\
\sum_{j_3=0}^{N_3-1}\sum_{j_2=0}^{N_2-1}\sum_{j_1=0}^{N_1-1}c_{0,j_3}b_{1,j_2}a_{N_1-1,j_1}{U}_{j_1,j_2,j_3}\\
\vdots\\
\sum_{j_3=0}^{N_3-1}\sum_{j_2=0}^{N_2-1}\sum_{j_1=0}^{N_1-1}c_{1,j_3}b_{0,j_2}a_{0,j_1}{U}_{j_1,j_2,j_3}\\
\vdots\\
\sum_{j_3=0}^{N_3-1}\sum_{j_2=0}^{N_2-1}\sum_{j_1=0}^{N_1-1}c_{1,j_3}b_{0,j_2}a_{N_1-1,j_1}{U}_{j_1,j_2,j_3}\\
\vdots\\
\sum_{j_3=0}^{N_3-1}\sum_{j_2=0}^{N_2-1}\sum_{j_1=0}^{N_1-1}c_{N_3-1,j_3}b_{N_2-1,j_2}a_{0,j_1}{U}_{j_1,j_2,j_3}\\
\vdots\\
\sum_{j_3=0}^{N_3-1}\sum_{j_2=0}^{N_2-1}\sum_{j_1=0}^{N_1-1}c_{N_3-1,j_3}b_{N_2-1,j_2}a_{N_1-1,j_1}{U}_{j_1,j_2,j_3}\\
 \end{array}
 \right).
\end{align*}
This completes the proof.}
\end{prf}

Under the periodic boundary condition \eqref{DNLS-PBS}, then, according to Lemma \ref{DNLS:lem:2.1}, we have
\begin{align}\label{finite-difference-method}
|{\bm U}|_{1,h}^2&={\langle-\big({\bm I}_{N_{3}}\otimes {\bm I}_{N_{2}}\otimes{\bm B}_{1}+{\bm I}_{N_{3}}\otimes {\bm B}_{2}\otimes {\bm I}_{N_{1}}
+{\bm B}_{3}\otimes {\bm I}_{N_{2}}\otimes {\bm I}_{N_{1}}\big){\bm U},{\bm U}\rangle_{h}}\nonumber\\
&:={\langle-{\bm\Delta}_{1,h}{\bm U},{\bm U}\rangle_{h}},
\end{align}
where {${\bm I}_{N_r}$} is an identity matrix of order $N_r$, and
\begin{align*}
{\bm B}_r=\frac{1}{h_r^2}\left(\begin{array}{ccccc}
              -2& 1 &0 &\cdots &1 \\
              1& -2 &1& \cdots &0\\
              \vdots &\vdots &\ddots& \ddots &\vdots\\
              0&0& \cdots&-2&1\\
              1&  0& \cdots&1 & -2\\
             \end{array}
\right)_{N_r\times N_r},r=1,2,3.
\end{align*}
Here, ${\bm B}_r,\ r=1,2,3$ is the usual discretization of the second partial derivative, by taking into
account of the periodic boundary condition. 

%

\subsection{Fourier pseudo-spectral method and some Lemmas }
xWe define
\begin{align*}
&S_{N}^{'''}=\text{span}\{g_{j_1}(x)g_{j_2}(y)g_{j_3}(z),\ 0\leq j_r\leq N_r-1, r=1,2,3\}
\end{align*}
as the interpolation space, where $g_{j_1}(x)$, $g_{j_2}(y)$ and $g_{j_3}(z)$ are trigonometric polynomials of degree $N_{1}/2$, $N_{2}/2$ and $N_{3}/2$,
given respectively by
\begin{align*}
  &g_{j_1}(x)=\frac{1}{N_{1}}\sum_{l=-N_{1}/2}^{N_{1}/2}\frac{1}{a_{l}}e^{\text{i}l\mu_{x} (x-x_{j_1})},\ a_{l}=\left \{
 \aligned
 &1,\ |l|<\frac{N_1}{2},\\
 &2,\ |l|=\frac{N_1}{2},
 \endaligned
 \right.\ \mu_x=\frac{2\pi}{l_1},\\
 & g_{j_2}(y)=\frac{1}{N_{2}}\sum_{p=-N_{2}/2}^{N_{2}/2}\frac{1}{b_{p}}e^{\text{i}p\mu_{y} (y-y_{j_2})}, b_{p}=\left \{
 \aligned
 &1,\ |p|<\frac{N_2}{2},\\
 &2,\ |p|=\frac{N_2}{2},
 \endaligned
 \right.\ \mu_y=\frac{2\pi}{l_2},\\
  &g_{j_3}(z)=\frac{1}{N_{3}}\sum_{q=-N_{3}/2}^{N_{3}/2}\frac{1}{c_{q}}e^{\text{i}q\mu_{z} (z-z_{j_3})},\ c_{q}=\left \{
 \aligned
 &1,\ |q|<\frac{N_3}{2},\\
 &2,\ |q|=\frac{N_3}{2},
 \endaligned
 \right.\ \mu_z=\frac{2\pi}{l_3},
\end{align*}
 {where $l_1,\ l_2$ and $l_3$ are the lengths of the $x$-,\ $y$- and $z$-directions, respectively, for the computational domain $\Omega$.}

{Denoting $C(\Omega)$ as the space of continuous functions on $\Omega\subset\mathbb{R}^3$,} we then define the interpolation operator $I_{N}: C(\Omega)\to S_{N}^{'''}$ as follows:
\begin{align}\label{Interpolation-operator}
I_{N}U(t,x,y,z)=\sum_{j_1=0}^{N_{1}-1}\sum_{j_2=0}^{N_{2}-1}\sum_{j_3=0}^{N_{3}-1}U_{j_1,j_2,j_3}g_{j_1}(x)g_{j_2}(y)g_{j_3}(z),
\end{align}
where $U_{j_1,j_2,j_3}=U(t,x_{j_1},y_{j_2},z_{j_3}),\ g_{j_1}(x_j)=\delta_{j_1}^j,\ g_{j_2}(x_k)=\delta_{j_2}^k$ and $g_{j_3}(x_l)=\delta_{j_2}^l$. 
Taking the second partial derivative with respect to variable $x$, and then evaluating the resulting expressions at the collocation point ($x_{j_1},y_{j_2},z_{j_3}$), we obtain
{\begin{align}\label{DNLS-2.3}
\frac{\partial^{2} I_{N}U(t,x_{j_1},y_{j_2},z_{j_3})}{\partial x^{2}}&=\sum_{j=0}^{N_{1}-1}\sum_{k=0}^{N_{2}-1}\sum_{m=0}^{N_{3}-1}U_{j,k,m}\frac{d^{2}g_{j}(x_{j_1})}{dx^{2}}g_{k}(y_{j_2})g_{m}(z_{j_3})\nonumber\\
&=\sum_{j=0}^{N_{1}-1}U_{j,j_2,j_3}\frac{d^{2}g_{j}(x_{j_1})}{dx^{2}}\nonumber\\
&=\sum_{j=0}^{N_{1}-1}({\bm D}_{2}^{x})_{j_1,j}U_{j,j_2,j_3},
\end{align}
where ${\bm D}_{2}^{x}$ is a real symmetric matrix of order $N_1$ with elements given by
\begin{align*}
({\bm D}_{2}^{x})_{j_1,j}=\frac{d^{2}g_{j}(x_{j_1})}{dx^{2}}.
\end{align*}
Similarly, we can obtain
\begin{align}\label{DNLS-2.4}
\frac{\partial^{2} I_{N}U(t,x_{j_1},y_{j_2},z_{j_3})}{\partial y^{2}}&=\sum_{k=0}^{N_{2}-1}U_{j_1,k,j_3}\frac{d^{2}g_{k}(y_{j_2})}{dy^{2}}
=\sum_{k=0}^{N_{2}-1}U_{j_1,k,j_3}({\bm D}_2^y)_{j_2,k},\\\label{DNLS-2.5}
 \frac{\partial^{2} I_{N}U(t,x_{j_1},y_{j_2},z_{j_3})}{\partial z^{2}}&=\sum_{m=0}^{N_{3}-1}U_{j_1,j_2,m}\frac{d^{2}g_{m}(z_{j_3})}{dz^{2}}
=\sum_{m=0}^{N_{3}-1}U_{j_1,j_2,m}({\bm D}_2^z)_{j_3,m},
\end{align}}
where ${\bm D}_{2}^{y}$ and ${\bm D}_{2}^{z}$ are real symmetric matrices of order $N_r,\ r=2,3$, respectively, with elements given by
\begin{align*}
({\bm D}_{2}^{y})_{j_2,k}=\frac{d^{2}g_{k}(y_{j_2})}{dy^{2}},\ ({\bm D}_2^z)_{j_3,m}=\frac{d^{2}g_{m}(z_{j_3})}{dz^{2}}.
\end{align*}
{According to Lemma \ref{DNLS:lem:2.1} and \eqref{Interpolation-operator}, we can deduce from \eqref{DNLS-2.3}-\eqref{DNLS-2.5} that
\begin{align}\label{DNLS-2.6}
\left(\begin{array}{c}
 \Delta I_{N}U(t,x_{0},y_{0},z_{0})\\
              \Delta I_{N}U(t,x_{1},y_{0},z_{0}) \\
              \vdots \\
             \Delta I_{N}U(t,x_{N_1-1},y_{0},z_{0}) \\
             \vdots \\
              \Delta I_{N}U(t,x_{0},y_{N_2-1},z_{N_3-1})\\
              \vdots \\
              \Delta I_{N}U(t,x_{N_1-1},y_{N_2-1},z_{N_3-1})\\
\end{array}
\right)={\bm \Delta}_{h}{\bm U},\ {\bm U}\in\mathbb{V}_h,
\end{align}
where
\begin{align}\label{Laplace-discrete}
{\bm \Delta}_{h}= {\bm I}_{N_{3}}\otimes {\bm I}_{N_{2}}\otimes{\bm D}_{2}^{x}+{\bm I}_{N_{3}}\otimes {\bm D}_{2}^{y}\otimes {\bm I}_{N_{1}}+{\bm D}_{2}^{z}\otimes {\bm I}_{N_{2}}\otimes {\bm I}_{N_{1}}.
\end{align}}
 Then, we introduce a new semi-norm induced by ${\bm \Delta}_h$ for ${\bm U}\in\mathbb{V}_{h}$, as follows:
\begin{align*}
 |{\bm U}|_{h}^2={\langle-{\bm \Delta}_{h}{\bm U},{\bm U}\rangle_{h}}.
\end{align*}


{ \begin{lem}\rm \label{3d_DNLSlem2.1}\cite{GCW14,HNO06} For the matrices ${\bm B}_{r},r=1,2,3$ and ${\bm D}_{2}^{w},w=x,y,z$, the following results hold
\begin{align*}
&{\bm B}_{r}={\bm F}_{N_r}^{H}{\bm\Lambda}_{r}{\bm F}_{N_r},\ {\bm\Lambda}_{r}=-\frac{4}{h_{r}^{2}}\text{\rm diag}\Big[\sin^{2}\frac{0}{N_{r}}\pi,\sin^{2}\frac{1}{N_{r}}\pi,\cdots,\sin^{2}\frac{N_r-1}{N_{r}}\pi\Big],\\
&{\bm D}_{2}^{x}={\bm F}_{N_1}^{H}{\bm\Lambda}_{x}{\bm F}_{N_1},\ {\bm\Lambda}_{x}=-\mu_x^2\text{\rm diag}\Big[0^2,1^2,\cdots,\big(\ \frac{N_1}{2}\big)^2,\big(-\frac{N_1}{2}+1\big)^2,\cdots,(-2)^2,(-1)^2\Big],\\
&{\bm D}_{2}^{y}={\bm F}_{N_2}^{H}{\bm\Lambda}_{y}{\bm F}_{N_2},\ {\bm\Lambda}_{y}=-\mu_y^2\text{\rm diag}\Big[0^2,1^2,\cdots,\big(\ \frac{N_2}{2}\big)^2,\big(-\frac{N_2}{2}+1\big)^2,\cdots,(-2)^2,(-1)^2\Big],\\
& {\bm D}_{2}^{z}={\bm F}_{N_3}^{H}{\bm\Lambda}_{z}{\bm F}_{N_3},\ \ {\bm\Lambda}_{z}=-\mu_z^2\text{\rm diag}\Big[0^2,1^2,\cdots,\big(\ \frac{N_3}{2}\big)^2,\big(-\frac{N_3}{2}+1\big)^2,\cdots,(-2)^2,(-1)^2\Big],
\end{align*}
where  ${\bm F}_{N_r},\ r=1,2,3$ is the discrete Fourier matrix of order $N_r$.
In addition, according to Definition \ref{DNLS:den2.1}, Lemma \ref{DNLS:lem:2.1} and the inequality on the eigenvalues of ${\bm B}_r,\ r=1,2,3$ and ${\bm D}_2^w,w=x,y,z$ (for the proof of the inequalities, please refer to Ref. \cite{GWWC17}), we have
\begin{align}\label{3d_NLS:2.10}
&\frac{4}{\pi^{2}}\langle-\big({\bm I}_{N_{3}}\otimes {\bm I}_{N_{2}}\otimes{\bm\Lambda}_{x}\big){{\bm U}},{{\bm U}}\rangle_{h}\le\langle-\big({\bm I}_{N_{3}}\otimes {\bm I}_{N_{2}}\otimes{\bm\Lambda}_{1}\big){{\bm U}},{{\bm U}}\rangle_{h}\nonumber\\
&~~~~~~~~~~~~~~~~~~~~~~~~~~~~~~~~~~~~~~~~~~~~~~~~~~~~~~~\le\langle-\big({\bm I}_{N_{3}}\otimes {\bm I}_{N_{2}}\otimes{\bm\Lambda}_{x}\big){{\bm U}},{{\bm U}}\rangle_{h},\\\label{3d_NLS:2.11}
&\frac{4}{\pi^{2}}\langle-\big({\bm I}_{N_{3}}\otimes {\bm\Lambda}_{y}\otimes {\bm I}_{N_{1}}\big){{\bm U}},{{\bm U}}\rangle_{h}\le\langle-\big({\bm I}_{N_{3}}\otimes{\bm\Lambda}_{2}\otimes{\bm I}_{N_{1}}\big){{\bm U}},{{\bm U}}\rangle_{h}\nonumber\\
&~~~~~~~~~~~~~~~~~~~~~~~~~~~~~~~~~~~~~~~~~~~~~~~~~~~~~~~\le\langle-\big({\bm I}_{N_{3}}\otimes{\bm\Lambda}_{y}\otimes{\bm I}_{N_{1}}\big){{\bm U}},{{\bm U}}\rangle_{h},\\\label{3d_NLS:2.12}
&\frac{4}{\pi^{2}}\langle-\big({\bm\Lambda}_{z}\otimes {\bm I}_{N_{2}}\otimes{\bm I}_{N_{1}}\big){{\bm U}},{{\bm U}}\rangle_{h}\le\langle-\big(\textcolor{blue}{\Lambda_{3}}\otimes {\bm I}_{N_{2}}\otimes{\bm I}_{N_{1}}\big){{\bm U}},{{\bm U}}\rangle_{h}\nonumber\\
&~~~~~~~~~~~~~~~~~~~~~~~~~~~~~~~~~~~~~~~~~~~~~~~~~~~~~~~\le\langle-\big({\bm\Lambda}_{z}\otimes {\bm I}_{N_{2}}\otimes{\bm I}_{N_{1}}\big){{\bm U}},{{\bm U}}\rangle_{h},
\end{align}
where ${\bm U}\in\mathbb{V}_h$.
\end{lem}}
Now, we give the following equivalence between $|{\bm U}|_{1,h}$ and $|{\bm U}|_{h}$.
\begin{lem}\rm \label{3d_DNLSlem2.2} For any grid function ${\bm U} \in \mathbb{V}_{h}$, we have
\begin{align*}
&|{\bm U}|_{1,h}\leq |{\bm U}|_{h}\leq\frac{\pi}{2}|{\bm U}|_{1,h}.
\end{align*}
\end{lem}
\begin{prf}\rm
Denoting by
\begin{align*}
I^{2}:=|{\bm U}|_{h}^{2}&=\langle-\Delta_{h}{\bm U},{\bm U}\rangle_{h}\nonumber\\
&=\langle-\big({\bm I}_{N_{3}}\otimes {\bm I}_{N_{2}}\otimes{\bm D}_{2}^{x}\big){\bm U},{\bm U}\rangle_{h}+\langle-\big({\bm I}_{N_{3}}\otimes {\bm D}_{2}^{y}\otimes {\bm I}_{N_{1}}\big){\bm U},{\bm U}\rangle_{h}\nonumber\\
&+\langle-\big({\bm D}_{2}^{z}\otimes {\bm I}_{N_{2}}\otimes {\bm I}_{N_{1}}\big){\bm U},{\bm U}\rangle_{h}:=I_{1}^{2}+I_{2}^{2}+I_{3}^{2},
\end{align*}
and
\begin{align*}
J^{2}:=|{\bm U}|_{1,h}^{2}&=\langle-\Delta_{1,h}{\bm U},{\bm U}\rangle_{h}\nonumber\\
&=\langle-\big({\bm I}_{N_{3}}\otimes {\bm I}_{N_{2}}\otimes{\bm B}_{1}\big){\bm U},{\bm U}\rangle_{h}+\langle-\big({\bm I}_{N_{3}}\otimes {\bm B}_{2}\otimes {\bm I}_{N_{1}}\big){\bm U},{\bm U}\rangle_{h}\nonumber\\
&+\langle-\big({\bm B}_{3}\otimes {\bm I}_{N_{2}}\otimes {\bm I}_{N_{1}}\big){\bm U},{\bm U}\rangle_{h}:=J_1^{2}+J_2^{2}+J_3^{2}.
\end{align*}
{With Lemmas \ref{DNLS:lem:2.1}-\ref{3d_DNLSlem2.1} and Corollary \ref{DNLS:cor:2.1}, we obtain
\begin{align*}
I_{1}^{2}&=\langle-\big({\bm I}_{N_{3}}\otimes {\bm I}_{N_{2}}\otimes{\bm F}_{N_1}^{H}{\bm\Lambda}_{x}{\bm F}_{N_1}\big){\bm U},{\bm U}\rangle_{h}\\
&=\langle-\big({\bm F}_{N_3}^{H}{\bm F}_{N_3}\otimes {\bm F}_{N_2}^{H}{\bm F}_{N_2}\otimes{\bm F}_{N_1}^{H}{\bm\Lambda}_{x}{\bm F}_{N_1}\big){\bm U},{\bm U}\rangle_{h}\\
&=\langle-\big({\bm F}_{N_3}\otimes {\bm F}_{N_2}\otimes{\bm F}_{N_1}\big)^{H}\big({\bm F}_{N_3}\otimes {\bm F}_{N_2}\otimes{\bm\Lambda}_{x}{\bm F}_{N_1}\big){\bm U},{\bm U}\rangle_{h}\\
&=\langle-\big({\bm F}_{N_3}\otimes {\bm F}_{N_2}\otimes{\bm\Lambda}_{x}{\bm F}_{N_1}\big){\bm U},\big({\bm F}_{N_3}\otimes {\bm F}_{N_2}\otimes{\bm F}_{N_1}\big){\bm U}\rangle_{h}\\
&=\langle-\big({\bm I}_{N_{3}}\otimes {\bm I}_{N_{2}}\otimes{\bm\Lambda}_{x}\big)\big({\bm F}_{N_3}\otimes {\bm F}_{N_2}\otimes{\bm F}_{N_1}\big){{\bm U}},\big({\bm F}_{N_3}\otimes {\bm F}_{N_2}\otimes{\bm F}_{N_1}\big){{\bm U}}\rangle_{h}\\
&=\langle-\big({\bm I}_{N_{3}}\otimes {\bm I}_{N_{2}}\otimes{\bm\Lambda}_{x}\big)\widetilde{{\bm U}},\widetilde{{\bm U}}\rangle_{h},
\end{align*}
where $\widetilde{{\bm U}}=\big({\bm F}_{N_3}\otimes {\bm F}_{N_2}\otimes{\bm F}_{N_1}\big){\bm U}$.} By the similar argument, we have
\begin{align*}
&I_{2}^{2}=\langle-\big({\bm I}_{N_{3}}\otimes {\bm\Lambda}_{y}\otimes{\bm I}_{N_{1}}\big)\widetilde{{\bm U}},\widetilde{{\bm U}}\rangle_{h},\ I_{3}^{2}=\langle-\big({\bm\Lambda}_{z}\otimes {\bm I}_{N_{2}}\otimes{\bm I}_{N_{1}}\big)\widetilde{{\bm U}},\widetilde{{\bm U}}\rangle_{h},\\
&J_1^{2}=\langle-\big({\bm I}_{N_{3}}\otimes {\bm I}_{N_{2}}\otimes{\bm\Lambda}_{1}\big)\widetilde{{\bm U}},\widetilde{{\bm U}}\rangle_{h},\ J_2^{2}=\langle-\big({\bm I}_{N_{3}}\otimes {\bm\Lambda}_{2}\otimes{\bm I}_{N_{1}}\big)\widetilde{{\bm U}},\widetilde{{\bm U}}\rangle_{h},\\ &J_3^{2}=\langle-\big({\bm\Lambda}_{3}\otimes {\bm I}_{N_{2}}\otimes{\bm I}_{N_{1}}\big)\widetilde{{\bm U}},\widetilde{{\bm U}}\rangle_{h}.
\end{align*}
With the use of \eqref{3d_NLS:2.10}-\eqref{3d_NLS:2.12}, we can get
\begin{align*}
\frac{4}{\pi^2}I_{r}^{2}\leq J_r^{2}\leq I_r^{2},\  r=1,2,3,
\end{align*}
which implies that
\begin{align}\label{3d_NLS:2.24}
J_{r}^{2}\leq I_r^{2}\leq \frac{\pi^2}{4}J_r^{2}, \ r=1,2,3.
\end{align}
Then, we deduce from \eqref{3d_NLS:2.24} that
\begin{align*}
J^{2}\leq I^{2}\leq \frac{\pi^2}{4}J^{2},
\end{align*}
that is,
\begin{align*}
|{\bm U}|_{1,h}\leq |{\bm U}|_{h}\leq\frac{\pi}{2}|{\bm U}|_{1,h}.
\end{align*}
This completes the proof.
\qed
\end{prf}

\begin{lem}\rm \cite{zhou90} \label{3d_DNLSlem2.3} For any grid function ${\bm U}\in \mathbb{V}_{h}$, there are
\begin{align*}
{ \|{\bm U}\|_{h,p}\leq C\|{\bm U}\|_{h}^{\frac{3}{p}-\frac{1}{2}}\{|{\bm U}|_{1,h}+\|{\bm U}\|_{h}\}^{\frac{3}{2}-\frac{3}{p}}},
\end{align*}
for $2\leq p\leq 6$. And
\begin{align*}
\|{\bm U}\|_{h,p}\leq Ch^{\frac{3}{p}-\frac{3}{2}}\|{\bm U}\|_{h},
\end{align*}
for $2\leq p\leq \infty$, where $C$ is a constant independent of $h$ and the grid function ${\bm U}$.
\end{lem}
\begin{lem}\rm \label{3d_DNLSlem2.4} For any grid function ${\bm U}^{n}\in\mathbb{V}_h$, it holds
\begin{align*}
\|{\bm U}^{m}\|_{h}\leq 2\sum_{k=1}^{m-1}\|\widehat{\bm U}^{k}\|_{h}+\|{\bm U}^{1}\|_{h}+\|{\bm U}^{0}\|_{h},\ 1\leq m\leq M.
\end{align*}
\end{lem}
\begin{prf}\rm
By using the triangular inequality, we have
\begin{align*}
\|{\bm U}^{m}\|_{h}-\|{\bm U}^{m-2}\|_{h}\leq 2\|\widehat{\bm U}^{m-1}\|_{h}.
\end{align*}
Summing up for $m$ from 2 to $K$ and then replacing $K$ by $m$, we can obtain that
\begin{align*}
\|{\bm U}^{m}\|_{h}+\|{\bm U}^{m-1}\|_{h}\leq 2\sum_{k=2}^{m}\|\widehat{\bm U}^{k-1}\|_{h}+\|{\bm U}^{1}\|_{h}+\|{\bm U}^{0}\|_{h},
\end{align*}
which further implies that
\begin{align*}
\|{\bm U}^{m}\|_{h}&\leq 2\sum_{k=1}^{m-1}\|\widehat{\bm U}^{k}\|_{h}+\|{\bm U}^{1}\|_{h}+\|{\bm U}^{0}\|_{h}.
\end{align*}

\end{prf}

 For brevity, we denote $U_{\vec{j}}^n$, $\Psi_{\vec{j}}^n$, $u_{\vec{j}}^n$ and $\psi_{\vec{j}}^n$ as the numerical approximations and the exact
solutions of $u(t,x,y,z)$ and $\psi(t,x,y,z)$ at the grid point $(t_n,x_{j_1},y_{j_2},z_{j_3})$, respectively. Note that $\Psi_{\vec{j}}^n=e^{-\gamma t_n}U_{\vec{j}}^n$ and $\psi_{\vec{j}}^n=e^{-\gamma t_n}u_{\vec{j}}^n$. Throughout the paper, let $C$ be a generic positive constant independent of $h$ and $\tau$, which may be different in different case.

\subsection{A linearly implicit and conservative scheme}\label{3D_NLS:Sec.2.1}
Discretizing the system \eqref{3DNLD-eq-1.2} using in time a linearly implicit structure-preserving method \cite{ABB13} and in space
the standard Fourier pseudo-spectral method, we obtain a linearly implicit and conservative Fourier pseudo-spectral (denoted by LI-CFP) scheme, as follows:
\begin{align}\label{3d_NLS:eq2.3}
&\text{i}\delta_{t}{\bm U}^{n}+ {\bm \Delta}_{h}\widehat{\bm U}^{n}+\beta e^{-2\gamma t_{n}}\big|{\bm U}^{n}\big|^{2}\cdot{ \widehat{\bm U}}^{n}=0,\ {\bm U}^{n}\in \mathbb{V}_{h},\ n=1,\cdots,M-1,
\end{align}
where $|{\bm U}^{n}|^{2}={\bm U}^{n}\cdot \bar{{\bm U}}^{n}$. Since the scheme \eqref{3d_NLS:eq2.3} is a three-level, {${\bm U}^{1}$ is obtained by the following modified Crank-Nicolson scheme
\begin{align}\label{3d_NLS:eq2.4}
&\text{i}\delta_{t}^+{\bm U}^{0}+{\bm \Delta}_{h}{\bm U}^{\frac{1}{2}}+\beta e^{-2\gamma t_{\frac{1}{2}}}\big|{\bm U}^{(\frac{1}{2})}\big|^{2}\cdot{\bm U}^{\frac{1}{2}}=0,\\\label{3d_NLS-eq2.5}
&{\bm U}_{\vec{j}}^{(\frac{1}{2})}=u(0,{\bm x}_j)+\frac{\tau}{2}\text{i}\Big({\Delta}u(0,{\bm x}_j)+\beta e^{-2\gamma t_{0}} |u(0,{\bm x}_j)|^2u(0,{\bm x}_j)\Big),\ \vec{j}\in{J}_h,
\end{align}
where $\Delta=\partial_{xx}+\partial_{yy}+\partial_{zz}$.}

\subsection{Conservation properties of the scheme}

In this subsection, we show that the proposed scheme can preserve the discrete total mass and energy conservation laws, respectively.
\begin{thm}\rm \label{3d_DNLSthm2.1} The scheme \eqref{3d_NLS:eq2.3}-\eqref{3d_NLS:eq2.4} possesses the following discrete total mass conservation law
\begin{align}\label{mass-conservation-law}
&\mathcal{M}^{n}=\mathcal{M}^{0},\ \mathcal{M}^{n}:=\|{\bm U}^{n}\|_{h}^{2},\ n=1,2,\cdots,M.
\end{align}
\end{thm}
\begin{prf}\rm
 We make the discrete inner product of \eqref{3d_NLS:eq2.3} with $2\widehat{\bm U}^{n}$, and take the imaginary part of the resulting equation to arrive at
\begin{align}\label{3d_NLS:2.34}
\frac{1}{2\tau}\Big(\|{\bm U}^{n+1}\|_{h}^{2}-\|{\bm U}^{n-1}\|_{h}^{2}\Big)-&2\text{Im}\langle-\Delta_{h}{ \widehat{\bm U}}^{n}, \widehat{\bm U}^{n}\rangle_{h}\nonumber\\
&+2\beta e^{-2\gamma t_{n}}\text{Im}\langle |{\bm U}^{n}|^2\cdot\widehat{\bm U}^{n}, \widehat{\bm U}^{n}\rangle_{h}=0.
\end{align}
Thanks to
\begin{align*}
&\text{Im}\langle-\Delta_{h}\widehat{\bm U}^{n}, \widehat{\bm U}^{n}\rangle_{h}=0,\ \ \text{Im}\langle |{\bm U}^{n}|^2\cdot\widehat{\bm U}^{n}, {\widehat{\bm U}}^{n}\rangle_{h}=0,
\end{align*}
we can deduce from (\ref{3d_NLS:2.34}) that
\begin{align*}
\|{\bm U}^{n+1}\|_{h}^{2}=\|{\bm U}^{n-1}\|_{h}^{2},\ 1\leq n\leq M-1.
\end{align*}
An argument similar to \eqref{3d_NLS:eq2.4} used in \eqref{3d_NLS:2.34} shows that
\begin{align*}
\|{\bm U}^{1}\|_{h}^{2}=\|{\bm U}^{0}\|_{h}^{2}.
\end{align*}
This completes the proof.
\qed
\end{prf}


\begin{thm}\rm \label{3d_DNLSthm2.2}{The scheme \eqref{3d_NLS:eq2.3}} possesses the following discrete total energy conservation law
\begin{align}\label{LI-CFP-energy-conservation-law}
\mathcal{E}^{n}=\mathcal{E}^{0},\ n=1,\cdots,M-1,
\end{align}
where
\begin{align*}
\mathcal{E}^{n}=&\frac{1}{2}|{\bm U}^{n+1}|_{h}^{2}+\frac{1}{2}|{\bm U}^{n}|_{h}^{2}-\frac{\beta}{2}e^{-2\gamma t_{n}}h_{\Delta}\sum_{\vec{j}\in J_{h}}|U_{\vec{j}}^{n}|^{2}|U_{\vec{j}}^{n+1}|^{2}\nonumber\\
&-\frac{\beta}{2}\sum_{l=1}^{n}e^{-2\gamma t_{l-1}}(1-e^{-2\gamma\tau})h_{\Delta}\sum_{\vec{j}\in J_{h}}|U_{\vec{j}}^{l-1}|^{2}|U_{\vec{j}}^{l}|^{2}.
\end{align*}
\end{thm}
\begin{prf}\rm
We make the discrete inner product of \eqref{3d_NLS:eq2.3} with $-2\tau\delta_{t}{\bm U}^{n}$ and take the real part of the resulting equation to arrive at
\begin{align*}
\frac{1}{2}|{\bm U}^{n+1}|_{h}^{2}-\frac{1}{2}|{\bm U}^{n-1}|_{h}^{2}
-\frac{\beta}{2}e^{-2\gamma t_{n}}h_{\Delta}\sum_{\vec{j}\in J_{h}}\Big(|U_{\vec{j}}^{n}|^{2}|U_{\vec{j}}^{n+1}|^{2}
-|U_{\vec{j}}^{n}|^{2}|U_{\vec{j}}^{n-1}|^{2}\Big)=0.
\end{align*}
Summing up for $n$ from 1 to $m$ and then replacing $m$ by $n$, we finish the proof.
\end{prf}
\section{Unique solvability}
In this section, we show that the scheme \eqref{3d_NLS:eq2.3}-\eqref{3d_NLS:eq2.4} is uniquely solvable. For a fixed $n$, Eq. \eqref{3d_NLS:eq2.3} can be rewritten as the following equivalent form
\begin{align*}
{\bm A}^{n}\widehat{{\bm  U}}^{n}={{\bm  U}}^{n-1},\ {\bm A}^n={\bm I}+{\bm S}^n,\ n=1,2,\cdots,M-1,
\end{align*}
where ${\bm I}$ is the identity
matrix of order $N_1N_2N_3$ and ${\bm S}^n=-\text{i}\tau({\bm\Delta}_{h}+\beta e^{-2\gamma t_{n}}\text{diag}(|{\bm U}^{n}|^2))$ is a skew-Hermitian matrix depending on the solution at the previous time step. In order to obtain the unique
solvability of the proposed scheme, we need to prove that the matrix ${\bm A}^n$ is invertible

If ${\bm A}^n{\bm x}={\bm 0}$, we have
\begin{align*}
0={\bm x}^{H}{\bm A}^n{\bm x}={\bm x}^{H}{\bm x},
\end{align*}
where the above equality follows from the skew-Hermitian property of ${\bm S}^n$. Thus, ${\bm x}=0$, that is, ${\bm A}^n{\bm x}={\bm 0}$ has only zero solution. Therefore, ${\bm A}^{n}$ is invertible. An argument similar to \eqref{3d_NLS:eq2.4} as used above shows that ${\bm U}^1$ is uniquely solvable. For brevity, we omit the details.
\section{An a priori estimate}\label{Sec.4}
In this section, we will establish an optimal error estimate for the scheme \eqref{3d_NLS:eq2.3}-\eqref{3d_NLS:eq2.4} in discrete $L^{2}$-norm.
For simplicity, we let $\Omega=[0,2\pi]^{3}$. More general cuboid domain can be translated into $\Omega$. Let $C_{p}^{\infty}(\Omega)$ be a set of infinitely differentiable functions with the period $2\pi$ defined on $\Omega$ for all
variables and $H_{p}^{r}(\Omega)$ is the closure of $C_{p}^{\infty}(\Omega)$ in $H^{r}(\Omega)$.
The semi-norm and the norm of $H_{p}^{r}(\Omega)$ are denoted by $\arrowvert\cdot\arrowvert_{r}$
and $\|\cdot\|_{r}$ respectively.
$\|\cdot\|_{0}$ is denoted by $\|\cdot\|$ for simplicity. {We should note that $H_{p}^{0}(\Omega)$ is denoted by $L^2(\Omega)$, in this paper.}

Let $N_1=N_2=N_3=N$ (that is, $h_1=h_2=h_3=h$), the interpolation space $S_{N}^{'''}$ can be rewritten as
\begin{align*}
&S_{N}^{'''}=\Big\{ u| u=\sum_{\arrowvert j_1\arrowvert,\arrowvert j_2\arrowvert,\arrowvert j_3\arrowvert\leq\frac{N}{2}}
\frac{{\widehat u}_{j_1,j_2,j_3}}{c_{j_1}c_{j_2}c_{j_3}}e^{\text{i}(j_1x+j_2y+j_3z)}: {\widehat u}_{\frac{N}{2},j_2,j_3}={\widehat u}_{-\frac{N}{2},j_2,j_3},\\
&~~~~~~~~~~~~~~~~~~~{\widehat u}_{j_1,\frac{N}{2},j_3}={\widehat u}_{j_1,-\frac{N}{2},j_3},\
{\widehat u}_{j_1,j_2,\frac{N}{2}}={\widehat u}_{j_1,j_2,-\frac{N}{2}}\Big\},
\end{align*}
where $c_{l}=1,\ |l|<\frac{N}{2},\ c_{-\frac{N}{2}}=c_{\frac{N}{2}}=2$. The projection space is defined as
\begin{align*}
S_{N}=\Big\{u| u=\sum_{\arrowvert j_1\arrowvert,\arrowvert j_2\arrowvert,\arrowvert j_3\arrowvert\leq\frac{N}{2}}
{\widehat u}_{j_1,j_2,j_3}e^{\text{i}(j_1x+j_2y+j_3z)}\Big\}.
\end{align*}
It is clear that $S_{N-2}\subseteq S_{N}^{'''}\subseteq S_{N}$. Let $P_{N}:L^{2}(\Omega)\to S_{N}$ as the orthogonal projection operator,
and recall the interpolation operator $I_{N}:C(\Omega)\to S_{N}^{'''}$ {(see \eqref{Interpolation-operator})}.
Further, $P_{N}$ and $I_{N}$ satisfy \cite{GWWC17}:
\begin{align*}
&1.\  P_{N}\partial_{w}u=\partial_{w} P_{N}u,\  I_{N}\partial_{w}u\ne \partial_{w} I_{N}u,\ w=x,y,\ \text{or}\ z. \\
&2.\  P_{N}u=u,\ \forall u\in S_{N}, \ I_{N}u=u,\ \forall u\in S_{N}^{'''}.
\end{align*}
\begin{lem}\rm \label{3d_DNLSlem4.1} \cite{CWG16}
 For ${u}\in S_{N}^{'''}$,
$\Arrowvert { u}\Arrowvert\leq \Arrowvert {\bm u}\Arrowvert_{h}\leq 2\sqrt{2}\Arrowvert { u}\Arrowvert$.
\end{lem}
\begin{lem}\rm \label{3d_DNLSlem4.2}
 \cite{CQ82} If $0\leq l \leq s$ and ${u}\in H_{p}^{s}(\Omega)$, then
\begin{align*}
&\Arrowvert P_{N}{u}-{u}\Arrowvert_{l}\leq CN^{l-s}\arrowvert {u}\arrowvert_{s},\\
&||P_{N} u||_{l}\le C ||u||_{l},
\end{align*}
 and in addition if $s>3/2$ then
\begin{align*}
&\Arrowvert I_{N}{ u}-{ u}\Arrowvert_{l}\leq CN^{l-s}\arrowvert {u}\arrowvert_{s}.
\end{align*}
\end{lem}
\begin{lem}\rm \label{3d_DNLSlem4.3} \cite{JCWL17} For ${ u}\in H_{p}^{s}(\Omega),\ s>\frac{3}{2}$, let ${u}^{*}=P_{N-2}{u}$. Then, we have
 \begin{align*}
 \Arrowvert {\bm u}^{*}-{\bm u}\Arrowvert_{h}\leq CN^{-s} \arrowvert {u}\arrowvert_{s}.
 \end{align*}
\end{lem}

We rewrite \eqref{3DNLD-eq-1.2} as
\begin{align}\label{3DNLD_eq:1.2}
 \text{i}u_{t}+\Delta u+\beta e^{-2\gamma t}\big|u\big|^{2}u=0,
\end{align}
Let us denote by
\begin{align*}
u^{*}=P_{N-2}u,\ f^{*}=P_{N-2}f(u),
\end{align*}
where $f(u)=e^{-2\gamma t}\beta|u|^{2}u$. The projected equation of \eqref{3DNLD_eq:1.2} is
\begin{align}\label{3d_DNLSeq:4.9}
\text{i}\partial_{t} u^{*}+ \Delta  u^{*}+f^{*}=0.
\end{align}
We define
\begin{align}\label{3d_DNLSeq:4.10}
&\xi_{\vec{j}}^{n}=\text{i}\delta_{t}(u^{*})_{\vec{j}}^{n}+ \Delta(\widehat{{ u}^{*}})_{\vec{j}}^{n}
+\widehat{({f}^{*})_{\vec{j}}}^{n},\ n=1,\cdots,M-1,\\\label{3d_DNLSeq-4.4}
&\xi_{\vec{j}}^{0}=\text{i}\delta_{t}^+(u^{*})_{\vec{j}}^{0}+ \Delta({{ u}^{*}})_{\vec{j}}^{\frac{1}{2}}
+({f}^{*})^{\frac{1}{2}}_{\vec{j}},\ \vec{j}\in J_{h},
\end{align}
where $(u^{*})_{\vec{j}}^{n}=u^{*}(t_{n},{\bm x}_{\vec{j}})$ and $({ f}^{*})_{\vec{j}}^{n}=f^{*}(t_{n},{\bm x}_{\vec{j}})$.
By noting that $u(t,{\bm x})^{*}\in S_{N}^{'''}$, we obtain \begin{align}\label{3d_DNLSeq:4.11}
\Delta u^{*}(t_{n},{\bm x}_{\vec{j}})=\Delta (I_{N}u^{*}(t_{n},{\bm x}_{\vec{j}}))&=({\bm\Delta}_{h}({\bm u}^{*})^{n})_{\vec{j}},\ ({\bm u}^{*})^n\in\mathbb{V}_h,\ \vec{j}\in J_{h},
\end{align}
where
\begin{align*}
({\bm\Delta}_{h}({\bm u}^{*})^{n})_{\vec{j}}=\sum_{j=0}^{N_{1}-1}({\bm D}_{2}^{x})_{j_1,j}u_{j,j_2,j_3}^n+\sum_{k=0}^{N_{2}-1}u_{j_1,k,j_3}^n({\bm D}_2^y)_{j_2,k}+\sum_{m=0}^{N_{3}-1}u_{j_1,j_2,m}^n({\bm D}_2^z)_{j_3,m}.
\end{align*}. 
With \eqref{3d_DNLSeq:4.11}, we can deduce from \eqref{3d_DNLSeq:4.9}-\eqref{3d_DNLSeq-4.4}
\begin{align}
&\xi_{\vec{j}}^{n}=\text{i}\delta_{t}(u^{*})_{\vec{j}}^{n}-\text{i}\widehat{(\partial_{t}(u^{*}))}_{\vec{j}}^{n},\ n=1,\cdots,M-1,\\
&\xi_{\vec{j}}^{0}=\text{i}\delta_{t}^+(u^{*})_{\vec{j}}^{0}-\text{i}{(\partial_{t}(u^{*}))}_{\vec{j}}^{\frac{1}{2}}, \ \vec{j}\in J_{h}.
\end{align}
By the Taylor formula, we have
\begin{align}\label{3d_DNLS:4.14}
|\xi_{\vec{j}}^{n}|\le C\tau^2,\ n=0,1,\cdots,M.
\end{align}
We then define the error function by
\begin{align*}
{\eta}_{\vec{j}}^{n}=(u^{*})_{\vec{j}}^{n}-U_{\vec{j}}^{n},\ \vec{j}\in J_{h},\ n=0,1,\cdots,M.
\end{align*}
Subtracting \eqref{3d_NLS:eq2.3} and \eqref{3d_NLS:eq2.4} from \eqref{3d_DNLSeq:4.10} and \eqref{3d_DNLSeq-4.4}, respectively, we obtain the error equation
\begin{align}\label{3d_DNLS:4.16}
&\xi_{\vec{j}}^{n}=\text{i}\delta_{t}{\eta}_{\vec{j}}^{n}+( {\bm\Delta}_{h}\widehat{\bm \eta}^{n})_{\vec{j}}+G_{\vec{j}}^{n},\ n=1,\cdots,M-1,\\\label{3d_DNLS-4.16}
&\xi_{\vec{j}}^{0}=\text{i}\delta_{t}^+{\eta}_{\vec{j}}^{0}+( {\bm\Delta}_{h}{\bm \eta}^{\frac{1}{2}})_{\vec{j}}+G_{\vec{j}}^{0},\ \vec{j}\in J_{h},
\end{align}
where ${\bm \eta}^n\in\mathbb{V}_h$ and
\begin{align*}
G_{\vec{j}}^{0}={(f^{*})^{\frac{1}{2}}_{\vec{j}}}-\beta e^{-2\gamma t_{\frac{1}{2}}}|{U}_{\vec{j}}^{(\frac{1}{2})}|^{2}
{U}_{\vec{j}}^{\frac{1}{2}},\ G_{\vec{j}}^{n}=&\widehat{(f^{*})_{\vec{j}}}^{n}-\beta e^{-2\gamma t_{n}}|{U}_{\vec{j}}^{n}|^{2}
\widehat{U}_{\vec{j}}^{n}.
\end{align*}
For convenience, we rewrite \eqref{3d_DNLS:4.16} and \eqref{3d_DNLS-4.16} as
\begin{align}\label{3d_DNLS:4.17}
&{\bm \xi}^{n}=\text{i}\delta_{t}{\bm \eta}^{n}+{\bm\Delta}_{h}\widehat{\bm \eta}^{n}+{\bm G}^{n},\ n=1,\cdots,M-1,\\\label{3d_DNLS-error-equaion-2}
&{\bm\xi}^{0}=\text{i}\delta_{t}^+{\bm \eta}^{0}+{\bm\Delta}_{h}{\bm \eta}^{\frac{1}{2}}+{\bm G}^{0}.
\end{align}

\begin{lem}\rm \label{3d_DNLSlem4.5} We assume $u(t,{\bm x})\in C^{3}\Big(0,T;H_{p}^{s}(\Omega)\Big),\ s>\frac{3}{2}$. Then, there exists a constant $\tau_0>0$ sufficiently small,
 such that, when $0< \tau\leq \tau_0$, we have
\begin{align*}
&\|{\bm u}^{1}-{\bm U}^{1}\|_{h}\leq C(N^{-s}+\tau^{2}),\ \|{\bm \psi}^{1}-{\bm \Psi}^{1}\|_{h}\leq Ce^{-\gamma t_1}(N^{-s}+\tau^{2}),\\
\end{align*}
and
\begin{align*}
\tau\|{\bm \eta}^{\frac{1}{2}}\|_{1,h}^2\leq  C(N^{-s}+\tau^{2})^2.
\end{align*}
\end{lem}
\begin{prf}\rm
Denoting
\begin{align*}
&({\bm G}_{1})_{\vec{j}}^{0}={({ f}^{*})}_{\vec{j}}^{\frac{1}{2}}-{ f}_{\vec{j}}^{{\frac{1}{2}}},\ ({\bm G}_{2})_{\vec{j}}^{0}={ f}_{\vec{j}}^{\frac{1}{2}}-F(u (\frac{\tau}{2},{\bm x}_j),u_{\vec{j}}^{\frac{1}{2}}),\\
 &({\bm G}_{3})_{\vec{j}}^{0}=F({u}(\frac{\tau}{2},{\bm x}_j),u_{\vec{j}}^{\frac{1}{2}})-F({u}^{*}(\frac{\tau}{2},{\bm x}_j),({u}_{\vec{j}}^{*})^{\frac{1}{2}}),\\
& ({\bm G}_{4})_{\vec{j}}^{0}=F({u}^{*}(\frac{\tau}{2},{\bm x}_j),({u}_{\vec{j}}^{*})^{\frac{1}{2}})-F(U_{\vec{j}}^{(\frac{1}{2})},({u}^{*})_{\vec{j}}^{\frac{1}{2}}),\\
 &({\bm G}_{5})_{\vec{j}}^{0}=F(U_{\vec{j}}^{(\frac{1}{2})},({u}^{*})_{\vec{j}}^{\frac{1}{2}})-F(U_{\vec{j}}^{(\frac{1}{2})},{U}_{\vec{j}}^{\frac{1}{2}}),
\end{align*}
then, from \eqref{3d_DNLS:4.16}, we have
\begin{align}\label{3d_DNLS-4.13}
{\bm G}_{\vec{j}}^{0}=({\bm G}_{1})_{\vec{j}}^{0}+({\bm G}_{2})_{\vec{j}}^{0}+({\bm G}_{3})_{\vec{j}}^{0}+({\bm G}_{4})_{\vec{j}}^{0}+({\bm G}_{5})_{\vec{j}}^{0},
\end{align}
where
\begin{align*}
F(u,v)=\beta e^{-2\gamma t_{\frac{1}{2}}}|u|^{2}v,\ {\bm U}^{(\frac{1}{2})}=u(0,{\bm x}_j)+\text{i}\frac{\tau}{2}\Big({\Delta}u(0,{\bm x}_j)+|u(0,{\bm x}_j)|^2u(0,{\bm x}_j)\Big).
\end{align*}
With the triangle inequality, we have
\begin{align}\label{3d_DNLS-4.14}
\|{\bm G}^{0}\|_{h}\le\|{\bm G}_{1}^{0}\|_{h}+\|{\bm G}_{2}^{0}\|_{h}+\|{\bm G}_{3}^{0}\|_{h}+\|{\bm G}_{4}^{0}\|_{h}+\|{\bm G}_{5}^{0}\|_{h}.
\end{align}
With Lemma \ref{3d_DNLSlem4.3}, we have
\begin{align}\label{3d_DNLS-4.15}
\|{\bm G}_{1}^{0}\|_{h}\le CN^{-s},\ \|{\bm G}_{3}^{0}\|_{h}\le CN^{-s}.
\end{align}
By the Taylor formula, we obtain
\begin{align}\label{3d_DNLS-4.16}
\|{\bm G}_{2}^{0}\|_{h}\le C\tau^2.
\end{align}
Since
\begin{align}\label{3d_DNLS-4.17}
u(\frac{\tau}{2},{\bm x}_j)=u(0,{\bm x}_j)+\frac{\tau}{2}\partial_t(0,{\bm x}_j)+\frac{\tau^2}{4}\partial_{tt}(\frac{\theta\tau}{2},{\bm x}_j),\ 0<\theta<1,
\end{align}
where $\partial_t(0,{\bm x}_j)=\text{i}{\Delta}u(0,{\bm x}_j)+\text{i}|u(0,{\bm x}_j)|^2u(0,{\bm x}_j).$ Then, we have
\begin{align}\label{3d_DNLS-4.18}
|{\bm U}_{\vec{j}}^{(\frac{1}{2})}-u(\frac{\tau}{2},{\bm x}_j)|=\frac{\tau^2}{4}|\partial_{tt}(\frac{\theta\tau}{2},{\bm x}_j)|\le C\tau^2,\ 0<\theta<1,
\end{align}
which further implies that
\begin{align}\label{3d_DNLS-4.19}
|{ U}_{\vec{j}}^{(\frac{1}{2})}|\le|u(\frac{\tau}{2},{\bm x}_j)|+C\tau^2\le C,
\end{align}
that is
\begin{align}\label{3d_DNLS-L-infty}
\|{\bm U}^{(\frac{1}{2})}\|_{h,\infty}\le C.
\end{align}
Note that
\begin{align}\label{3d_DNLS-4.20}
\|{\bm U}^{(\frac{1}{2})}-{\bm u}^*(\frac{\tau}{2},{\bm x})\|_h&\le \|{\bm U}^{(\frac{1}{2})}-{\bm u}(\frac{\tau}{2},{\bm x})\|_h+\|{\bm u}(\frac{\tau}{2},{\bm x})-{\bm u}^*(\frac{\tau}{2},{\bm x})\|_h\nonumber\\
&\le C(N^{-s}+\tau^2),
\end{align}
we have
\begin{align}\label{3d_DNLS-4.21}
\|{\bm G}_{4}^{0}\|_{h}\le C(N^{-s}+\tau^2).
\end{align}
For ${\bm G}_5^0$, according to \eqref{3d_DNLS-L-infty}, we obtain
\begin{align}\label{3d_DNLS-4.22}
|({\bm G}_{5})_{\vec{j}}^{0}|&=|\beta e^{-2\gamma t_{\frac{1}{2}}}\big(|U_{\vec{j}}^{(\frac{1}{2})}|^{2}({u}^{*})_{\vec{j}}^{\frac{1}{2}}-|U_{\vec{j}}^{(\frac{1}{2})}|^{2}{U}_{\vec{j}}^{\frac{1}{2}}\big)|\nonumber\\
&\le C|({u}^{*})_{\vec{j}}^{\frac{1}{2}}-{U}_{\vec{j}}^{\frac{1}{2}}|\le C\big(|({u}^{*})_{\vec{j}}^{1}-{U}_{\vec{j}}^{1}|+|({u}^{*})_{\vec{j}}^{0}-{U}_{\vec{j}}^{0}|\big),
\end{align}
which further implies that
\begin{align}\label{3d_DNLS-4.23}
\|{\bm G}_{5}^0\|_h\le C\big(\|({\bm u}^{*})^{1}-{\bm U}^{1}\|_h+\|({\bm u}^{*})^{0}-{\bm U}^{0}\|_h\big)\le C\big(\|{\bm \eta}^{1}\|_h+N^{-s}\big).
\end{align}
Up to now, it follows from \eqref{3d_DNLS-4.15}-\eqref{3d_DNLS-4.16}, \eqref{3d_DNLS-4.21} and \eqref{3d_DNLS-4.23} that
\begin{align}\label{3d_DNLS-4.24}
\|{\bm G}^0\|_h\le C\|{\bm \eta}^{1}\|_h+C(N^{-s}+\tau^2).
\end{align}
We then make the inner product of \eqref{3d_DNLS-error-equaion-2} with $2{\bm \eta}^{\frac{1}{2}}$
\begin{align}\label{3d_DNLS-4.25}
\text{i}\delta_{t}^+\|{\bm \eta}^{0}\|_{h}^{2}-\frac{2}{\tau}\text{Im}\langle {\bm \eta}^{1},{\bm \eta}^{0}\rangle_{h}- 2|{\bm \eta}^{\frac{1}{2}}|_{h}^{2}+\langle {\bm G}^{0}, 2{\bm \eta}^{\frac{1}{2}}\rangle_{h}=\langle{\bm \xi}^{0},2{\bm \eta}^{\frac{1}{2}}\rangle_{h}.
\end{align}
The imaginary part of the above equation implies
\begin{align}\label{3d_DNLS-4.26}
\frac{1}{\tau}\big(\|{\bm \eta}^{1}\|_{h}^2-\|{\bm \eta}^{0}\|_{h}^{2}\big)+\text{Im}\langle {\bm G}^{0}, 2{\bm \eta}^{\frac{1}{2}}\rangle_{h}=\text{Im}\langle{\bm \xi}^{0},2{\bm \eta}^{\frac{1}{2}}\rangle_{h}.
\end{align}
 With the Cauchy-Schwarz inequality, triangular inequality, \eqref{3d_DNLS:4.14} and \eqref{3d_DNLS-4.24}, it is easy to see that
 \begin{align}\label{3d_DNLS-4.27}
 \|{\bm \eta}^{1}\|_{h}^{2}&\leq C\tau\|{\bm \eta}^{1}\|_{h}^{2}+C\|{\bm \eta}^{0}\|_{h}^{2}+C\tau\|{\bm \xi}^{0}\|_{h}^{2}+C\tau\|{\bm G}^{0}\|_{h}^{2}\nonumber\\
 &\le C(N^{-s}+\tau^2),
 \end{align}
 when $C\tau<1$.
  By virtue of Lemma \ref{3d_DNLSlem4.3} and equation \eqref{3d_DNLS-4.27}, we have
  \begin{align}\label{3d_DNLSeq:4.32}
 \|{\bm u}^{1}-{\bm U}^{1}\|_{h}^{2}\leq \|{\bm u}^{1}-({\bm u}^{*})^1\|_h^2+\|{\bm \eta}^{1}\|_{h}^{2}\leq C(N^{-s}+\tau^{2})^{2},
 \end{align}
which further implies that
\begin{align*}
\|{\bm u}^{1}-{\bm U}^{1}\|_{h}\leq C(N^{-s}+\tau^{2}),\ \|{\bm \psi}^{1}-{\bm \Psi}^{1}\|_{h}\leq e^{-\gamma t_1}C(N^{-s}+\tau^{2}).
\end{align*}

The real part of \eqref{3d_DNLS-4.25} reads
\begin{align*}
|{\bm \eta}^{\frac{1}{2}}|_{h}^{2}&=-\frac{1}{\tau}\text{Im}\langle {\bm \eta}^{1},{\bm \eta}^{0}\rangle_{h}+\frac{1}{2}\text{Re}\langle {\bm G}^{0}, 2{\bm \eta}^{\frac{1}{2}}\rangle_h-\frac{1}{2}\text{Re}\langle{\bm \xi}^{0},2{\bm \eta}^{\frac{1}{2}}\rangle_{h}.
\end{align*}
With \eqref{3d_DNLS-4.24} and the Cauchy-Schwarz inequality, we can get
\begin{align*}
\tau|{\bm \eta}^{\frac{1}{2}}|_{h}^{2}\leq C\big(\|{\bm \eta}^{1}\|_{h}^{2}+\|{\bm \eta}^{0}\|_{h}^{2}+\|{\bm G}^{0}\|_{h}^{2}+\|{\bm \xi}^{0}\|_{h}^{2}\big),
\end{align*}
which further implies that
 \begin{align*}
\tau|{\bm \eta}^{\frac{1}{2}}|_{h}^{2}\leq C(N^{-s}+\tau^{2})^{2},
 \end{align*}
 when $0< \tau\leq \tau_0$. According to Lemma \ref{3d_DNLSlem2.2}, we have
 \begin{align}\label{3d_DNLS-new-4.29}
\tau|{\bm \eta}^{\frac{1}{2}}|_{1,h}^{2}\leq C(N^{-s}+\tau^{2})^{2}.
 \end{align}
 This completes the proof.
\qed
\end{prf}
\begin{lem}\rm \label{3d_DNLSlem4.6} We assume $u(t,{\bm x})=C^{3}\Big(0,T;H_{p}^{s}(\Omega)\Big),\ s>\frac{3}{2}$.
Then, there exists a constant $\tau_0>0$ sufficiently small, such that when $0< \tau\leq \tau_0$, we have
\begin{align*}
\|{\bm u}^{2}-{\bm U}^{2}\|_{h}\leq C(N^{-s}+\tau^{2}),\ \|{\bm \psi}^{2}-{\bm \Psi}^{2}\|_{h}\leq Ce^{-\gamma t_2}(N^{-s}+\tau^{2}),
\end{align*}
and
\begin{align*}
\|{\bm \eta}^{2}\|_{h}^{2}+\tau\|\widehat{\bm \eta}^{1}\|_{1,h}^{2}\leq C(N^{-s}+\tau^{2})^{2}.
\end{align*}
\end{lem}
\begin{prf}\rm
Setting $n=1$ in \eqref{3d_DNLS:4.16}, we then multiply it by $2\widehat{\eta}_{\vec{j}}^{1}$ and sum them up for $\vec{j}\in J_{h}$ to obtain
\begin{align}\label{3d_DNLSeq-4.34}
\text{i}\delta_{t}\|{\bm \eta}^{1}\|_{h}^{2}-\frac{1}{\tau}\text{Im}\langle {\bm \eta}^{2},{\bm \eta}^{0}\rangle_{h}- 2|\widehat{\bm \eta}^{1}|_{h}^{2}+R^{1}=\langle{\bm \xi}^{1},2\widehat{\bm \eta}^{1}\rangle_{h},
\end{align}
where
\begin{align*}
R^{1}=2h_{\Delta}\sum_{\vec{j}\in J_{h}}\Big[\widehat{({f}^{*})}_{\vec{j}}^{1}-\beta e^{-2\gamma t_{1}}|{U}_{\vec{j}}^{1}|^{2}
\widehat{U}_{\vec{j}}^{1}\Big]\widehat{\bar{\eta}}_{\vec{j}}^{1}.
\end{align*}
$R^{1}$ can be written as
\begin{align}\label{3d_DNLS-neq-4.30}
R^{1}&=2\langle \widehat{{\bm f}^{*}}^{1}-\widehat{\bm f}^{1}, \widehat{\bm \eta}^{1}\rangle_{h}+ 2\langle \widehat{\bm f}^{1}-F({\bm u}^{1}), \widehat{\bm \eta}^{1}\rangle_{h}\nonumber\\
&+2\langle F({\bm u}^{1})-F(({\bm u}^{*})^{1}), \widehat{\bm \eta}^{1}\rangle_{h}+
2\langle F(({\bm u}^{*})^{1})-F({\bm U}^{1}), \widehat{\bm \eta}^{1}\rangle_{h}\nonumber\\
&:=R_{1}^{1}+R_{2}^{1}+R_{3}^{1}+R_{4}^{1},
\end{align}
where
\begin{align*}
F({\bm v}^1)=\beta e^{-2\gamma t_{1}}|{\bm v}^1|^{2}\widehat{\bm v}^1.
\end{align*}
Using the Taylor formula and Lemma \ref{3d_DNLSlem4.3}, we have
\begin{align}\label{3d_DNLSeq-4.44}
R_{1}^{1}\leq CN^{-2s}+\|\widehat{\bm \eta}^{1}\|_{h}^{2},\ R_{2}^{1}\leq C\tau^{4}+||\widehat{\bm \eta}^{1}||_{h}^{2}.
\end{align}
According to Lemma \ref{3d_DNLSlem4.2}, we have
\begin{align}\label{3d_DNLSeq-4.45}
\|({\bm u}^{*})^{n}\|_{h,\infty}&\leq \max\limits_{{\bm x}\in \Omega}|u^{*}(t_{n},{\bm x})|\leq C\|u^{*}(t_{n},{\bm x})\|_{2}\leq C\|u(t_{n},{\bm x})\|_{2}\leq C,
\end{align}
for $n=0,1,2$.
Then, with Lemma \ref{3d_DNLSlem4.3}, we get
\begin{align}\label{3d_DNLSeq-14.46}
R_{3}^{1}\leq CN^{-2s}+\|\widehat{\bm \eta}^{1}\|_{h}^{2}.
\end{align}
We rewrite $R_{4}^{1}$ by
\begin{align}\label{3d_DNLSeq-4.47}
R_{4}^{1}
&=2\beta e^{-2\gamma t_{1}} h_{\Delta}\sum_{\vec{j}\in J_{h}}\Big[|(u^{*})_{\vec{j}}^{1}|^{2}\widehat{(u^{*})}_{\vec{j}}^{1}
-|U_{\vec{j}}^{1}|^{2}\widehat{U}_{\vec{j}}^{1}\Big]\widehat{\bar{\eta}}_{\vec{j}}^{1}\nonumber\\
&=2\beta e^{-2\gamma t_{1}} h_{\Delta}\sum_{\vec{j}\in J_{h}}|(u^{*})_{\vec{j}}^{1}|^{2}|\widehat{\eta}_{\vec{j}}^{1}|^{2}\nonumber\\
&~~~+2\beta e^{-2\gamma t_{1}} h_{\Delta}\sum_{\vec{j}\in J_{h}}\Big(|(u^{*})_{\vec{j}}^{1}|^{2}-|U_{\vec{j}}^{1}|^{2}\Big)\Big(\widehat{(u^{*})}_{\vec{j}}^{1}
-\widehat{\eta}_{\vec{j}}^{1}\Big)\widehat{\bar{\eta}}_{\vec{j}}^{1}\nonumber\\
&:=R_{41}^{1}+R_{42}^{1}.
\end{align}
It is easy to see that
\begin{align}\label{3d_DNLSeq-4.48}
\text{Im}\big(R_{41}^{1}\big)=0,
\end{align}
and
\begin{align}\label{3d_DNLSeq-4.49}
|R_{42}^{1}|&\leq 2|\beta| e^{-2\gamma t_{1}} h_{\Delta}\sum_{\vec{j}\in J_{h}}\Big[|\eta_{\vec{j}}^{1}|^2+2|(u^{*})_{\vec{j}}^{1}|\cdot |\eta_{\vec{j}}^{1}|\Big]\Big(\Big|\widehat{(u^{*})}_{\vec{j}}^{1}\Big|
+\Big|\widehat{\eta}_{\vec{j}}^{1}\Big|\Big)\Big|\widehat{\eta}_{\vec{j}}^{1}\Big|\nonumber\\
&\le2|\beta| e^{-2\gamma t_{1}} h_{\Delta}\sum_{\vec{j}\in J_{h}}\Big[|\eta_{\vec{j}}^{1}|^2\Big|\widehat{\eta}_{\vec{j}}^{1}\Big|^2+
|\eta_{\vec{j}}^{1}|^2\Big|\widehat{(u^*)}_{\vec{j}}^{1}\Big|\Big|\widehat{\eta}_{\vec{j}}^{1}\Big|\nonumber\\
&+2\Big|(u^*)_{\vec{j}}^{1}\Big||\eta_{\vec{j}}^{1}|\Big|\widehat{\eta}_{\vec{j}}^{1}\Big|^2
+2\Big|(u^*)_{\vec{j}}^{1}\Big||\eta_{\vec{j}}^{1}|\Big|\widehat{(u^*)}_{\vec{j}}^{1}\Big|\Big|\widehat{\eta}_{\vec{j}}^{1}\Big|\nonumber\\
&\le Ch_{\Delta}\sum_{\vec{j}\in J_{h}}\Big[|\eta_{\vec{j}}^{1}|^2\Big|\widehat{\eta}_{\vec{j}}^{1}\Big|^2+|\eta_{\vec{j}}^1|^2+\Big|\widehat{\eta}_{\vec{j}}^{1}\Big|^2\Big]\nonumber\\
&\leq C\|{\bm \eta}^{1}\|_{h}^{2}+C\|\widehat{\bm \eta}^{1}\|_{h}^{2}+C\|{\bm \eta}^{1}\|_{h,3}^{2}\|\widehat{\bm \eta}^{1}\|_{h,6}^{2},
\end{align}
where the H\"{o}lder inequality is used.
With the above inequalities, we can prove that
\begin{align}\label{3d_DNLSeq-4.50}
|\text{Im}(R^1)|\leq C\|{\bm \eta}^{1}\|_{h}^{2}+C\|\widehat{\bm \eta}^{1}\|_{h}^{2}+C\|{\bm \eta}^{1}\|_{h,3}^{2}\|\widehat{\bm \eta}^{1}\|_{h,6}^{2}+C(N^{-s}+\tau^{2})^{2}.
\end{align}
The imaginary part of \eqref{3d_DNLSeq-4.34} reads
\begin{align}\label{3d_DNLSeq-4.51}
\frac{1}{2\tau}\Big(\|{\bm \eta}^{2}\|_{h}^{2}-\|{\bm \eta}^{0}\|_{h}^{2}\Big)+\text{Im}(R^{1})=\text{Im}\langle{\bm \xi}^{1},2\widehat{\bm \eta}^{1}\rangle_{h}.
\end{align}
Thanks to the Cauchy-Schwartz inequality, \eqref{3d_DNLS:4.14} and \eqref{3d_DNLSeq-4.50}, \eqref{3d_DNLSeq-4.51} reduces to
\begin{align}\label{3d_DNLSeq-4.52}
\frac{1}{\tau}\Big(\|{\bm \eta}^{2}\|_{h}^{2}-\|{\bm \eta}^{0}\|_{h}^{2}\Big)&\le C\|{\bm \eta}^{1}\|_{h}^{2}+C\|\widehat{\bm \eta}^{1}\|_{h}^{2}\nonumber\\
&+C\|{\bm \eta}^{1}\|_{h,3}^{2}\|\widehat{\bm \eta}^{1}\|_{h,6}^{2}+C(N^{-s}+\tau^{2})^{2}.
\end{align}
The real part of \eqref{3d_DNLSeq-4.34} reads
\begin{align}\label{3d_DNLSeq-4.53}
|\widehat{\bm \eta}^{1}|_{h}^{2}+\frac{1}{2\tau}\text{Im}\langle {\bm \eta}^{2},{\bm \eta}^{1}\rangle_{h}-\frac{1}{2}\text{Re}\big(R^{1}\big)=-\text{Re}\langle{\bm \xi}^{1},\widehat{\bm \eta}^{1}\rangle_{h}.
\end{align}
Combining  \eqref{3d_DNLSeq-4.45} with \eqref{3d_DNLSeq-4.47}, we obtain
\begin{align}\label{3d_DNLSeq-4.54}
|R_{41}^{1}|\leq C\|\widehat{\bm \eta}^{1}\|^{2}.
\end{align}
With Lemma \ref{3d_DNLSlem2.2}, the Cauchy-Schwartz inequality, \eqref{3d_DNLSeq-4.49} and \eqref{3d_DNLSeq-4.54}, \eqref{3d_DNLSeq-4.53} reduces to
\begin{align}\label{3d_DNLSeq-4.55}
|\widehat{\bm \eta}^{1}|_{1,h}^{2}&\leq C\|{\bm \eta}^{1}\|_{h}^{2}+C\|\widehat{\bm \eta}^{1}\|_{h}^{2}+C\|{\bm \eta}^{1}\|_{h,3}^{2}\|\widehat{\bm \eta}^{1}\|_{h,6}^{2}\nonumber\\
&+\frac{1}{4\tau}\Big(\|{\bm \eta}^{2}\|_{h}^{2}+\|{\bm \eta}^{0}\|_{h}^{2}\Big)+C(N^{-s}+\tau^{2})^{2}.
\end{align}

To reduce the nonlinear term in the above inequality, now we prove the following inequality
\begin{align}\label{3d_DNLSeq-4.56}
\|{\bm \eta}^{1}\|_{h,3}^{2}\|\widehat{\bm \eta}^{1}\|_{h,6}^{2}&\leq C\big(\|{\bm \eta}^{2}\|_{h}^{2}+\|{\bm \eta}^{1}\|_{h}^{2}+\|{\bm \eta}^{0}\|_{h}^{2} \big)+C\|\widehat{\bm \eta}^{1}\|_{h}^{2}+C(N^{-s}+\tau^{2})^{2},
\end{align}
with two different cases.

Firstly, we consider the case $\tau\leq h\ (h=\frac{2\pi}{N})$. We use Lemma \ref{3d_DNLSlem2.3} and \eqref{3d_DNLS-4.27} to get
\begin{align*}
&\|{\bm\eta}^{1}\|_{h,3}\leq Ch^{-\frac{1}{2}}\|{\bm \eta}^{1}\|_{h}\leq Ch^{\frac{3}{2}},\ \|\widehat{\bm\eta}^{1}\|_{6,h}\leq Ch^{-1}\|\widehat{\bm \eta}^{1}\|_{h}.
\end{align*}
When $ Ch\leq 1,$ we have
\begin{align}\label{3d_DNLSeq-4.59}
\|{\bm \eta}^{1}\|_{h,3}^{2}\|\widehat{\bm \eta}^{1}\|_{h,6}^{2}\leq C\|\widehat{\bm \eta}^{1}\|_{h}^{2}.
\end{align}

Secondly, we consider the case $\tau\ge h$. By Lemma \ref{3d_DNLSlem2.3}, we get
\begin{align}\label{3d_DNLSeq-4.60}
\|\widehat{\bm \eta}^{1}\|_{h,6}^{2}\leq C\big(|\widehat{\bm \eta}^{1}|_{1,h}+\|\widehat{\bm \eta}^{1}\|_{h})^{2}.
\end{align}
With \eqref{3d_DNLS-4.27} and \eqref{3d_DNLS-new-4.29}, together with Lemma \ref{3d_DNLSlem2.4} and the H\"older inequality, we have
\begin{align}\label{3d_DNLSeq-4.61}
\|{\bm \eta}^{1}\|_{h,3}^{2}&\leq \|{\bm \eta}^{1}\|_{h}\|{\bm \eta}^{1}\|_{h,6}\nonumber\\
&\leq \|{\bm \eta}^{1}\|_{h}\big(|{\bm \eta}^{1}|_{1,h}+\|{\bm \eta}^{1}\|_{h})\nonumber\\
&\leq \|{\bm \eta}^{1}\|_{h}\Big(2|{\bm \eta}^{\frac{1}{2}}|_{1,h}+|{\bm \eta}^{0}|_{1,h}+\|{\bm \eta}^{1}\|_{h}\Big)\nonumber\\
&\leq C\tau^{\frac{7}{2}}.
\end{align}

It follows from \eqref{3d_DNLSeq-4.60} and \eqref{3d_DNLSeq-4.61} that
\begin{align}\label{3d_DNLSeq-4.62}
\|{\bm \eta}^{1}\|_{h,3}^{2}\|\widehat{\bm \eta}^{1}\|_{h,6}^{2}&\leq \epsilon\tau\big(|\widehat{\bm \eta}^{1}|_{1,h}^{2}+\|\widehat{\bm \eta}^{1}\|_{h}^{2}\big),
\end{align}
when $C\tau^{\frac{5}{2}}\leq \epsilon$.
With the use of \eqref{3d_DNLSeq-4.62}, we can deduce from \eqref{3d_DNLSeq-4.55} that
\begin{align*}
|\widehat{\bm \eta}^{1}|_{1,h}^{2}&\leq C\|{\bm \eta}^{1}\|_{h}^{2}+C\|\widehat{\bm \eta}^{1}\|_{h}^{2}+\epsilon\tau|\widehat{\bm \eta}^{1}|_{1,h}^{2}\nonumber\\
&+\frac{1}{4\tau}\Big(\|{\bm \eta}^{2}\|_{h}^{2}+\|{\bm \eta}^{0}\|_{h}^{2}\Big)+C(N^{-s}+\tau^{2})^{2},
\end{align*}
which implies that
\begin{align}\label{3d_DNLSeq-4.63}
|\widehat{\bm \eta}^{1}|_{1,h}^{2}&\leq C\|{\bm \eta}^{1}\|_{h}^{2}+C\|\widehat{\bm \eta}^{1}\|_{h}^{2}+\frac{1}{4\tau}\Big(\|{\bm \eta}^{2}\|_{h}^{2}+\|{\bm \eta}^{0}\|_{h}^{2}\Big)+C(N^{-s}+\tau^{2})^{2},
\end{align}
when $\tau$ sufficiently small.
With \eqref{3d_DNLSeq-4.63}, \eqref{3d_DNLSeq-4.62} reduces to
\begin{align}\label{3d_DNLSeq-4.64}
\|{\bm \eta}^{1}\|_{h,3}^{2}\|\widehat{\bm \eta}^{1}\|_{h,6}^{2}&\leq C\big(\|{\bm \eta}^{2}\|_{h}^{2}+\|{\bm \eta}^{1}\|_{h}^{2}+\|{\bm \eta}^{0}\|_{h}^{2} \big)+C\|\widehat{\bm \eta}^{1}\|_{h}^{2}+C(N^{-s}+\tau^{2})^{2}.
\end{align}
Up to now, we have proved that \eqref{3d_DNLSeq-4.56} holds.
By using \eqref{3d_DNLSeq-4.56}, we can deduce from \eqref{3d_DNLSeq-4.52} that
\begin{align}\label{3d_DNLSeq-4.65}
\frac{1}{\tau}\Big(\|{\bm \eta}^{2}\|_{h}^{2}-\|{\bm \eta}^{0}\|_{h}^{2}\Big)&\le C\big(\|{\bm \eta}^{2}\|_{h}^{2}+\|{\bm \eta}^{1}\|_{h}^{2}+\|{\bm \eta}^{0}\|_{h}^{2}\big)+C(N^{-s}+\tau^{2})^{2}.
\end{align}
Then, we can prove
\begin{align}\label{3d_DNLSeq-4.66}
\|{\bm \eta}^{2}\|_{h}^{2}\le C(N^{-s}+\tau^{2})^{2},
\end{align}
where $\tau$ is sufficiently small, such that $C\tau\le \frac{1}{2}$. According to Lemma \ref{3d_DNLSlem4.3} and Eq. \eqref{3d_DNLSeq-4.66}, we have
  \begin{align}\label{3d_DNLSeq:4.32}
 \|{\bm u}^{2}-{\bm U}^{2}\|_{h}^{2}\leq \|{\bm u}^{2}-({\bm u}^{*})^2\|_h^2+\|{\bm \eta}^{2}\|_{h}^{2}\leq C(N^{-s}+\tau^{2})^{2},
 \end{align}
which further shows that
 \begin{align*}
  \|{\bm u}^{2}-{\bm U}^{2}\|_{h}\leq C(N^{-s}+\tau^{2}),\ \|{\bm \psi}^{2}-{\bm \Psi}^{2}\|_{h}\leq Ce^{-\gamma t_2}(N^{-s}+\tau^{2}).
 \end{align*}

Furthermore, from \eqref{3d_DNLSeq-4.55}-\eqref{3d_DNLSeq-4.56} and \eqref{3d_DNLSeq-4.66}, we see
\begin{align}\label{3d_DNLSeq-4.67}
\tau|\widehat{\bm \eta}^{1}|_{1,h}^{2}\leq C(N^{-s}+\tau^{2})^{2}.
\end{align}
This completes the proof.
\qed
\end{prf}

\begin{thm}\rm \label{3d_DNLSthm4.1}
We assume that the continuous solution $u(t,{\bm x})$ of (\ref{3DNLD_eq:1.2}) satisfies
\begin{align*}
u(t,{\bm x})\in C^{3}\Big(0,T;H_{p}^{s}(\Omega)\Big),\ s>\frac{3}{2}.
\end{align*}
Then, there exists a constant $s_0>0$, such that when $0<\tau,h\leq s_0$,
 we have the following error estimate for the scheme \eqref{3d_NLS:eq2.3}-\eqref{3d_NLS:eq2.4}
\begin{align*}
\|{\bm u}^{n}-{\bm U}^{n}\|_{h}\leq C_{0} (N^{-s}+\tau^{2}),\ 1\leq n\leq M,
\end{align*}
that is,
\begin{align*}
\|{\bm \psi}^{n}-{\bm \Psi}^{n}\|_{h}\leq C_{0} e^{-\gamma t_{n}}(N^{-s}+\tau^{2}),\ 1\leq n\leq M,
\end{align*}
where $C_{0}$ is a positive constant independent of $n,\tau$ and $h$.
\end{thm}
\subsection{The proof of Theorem 4.1}
In this section, we shall prove a slightly stronger result than Theorem 4.1, as follows:
\begin{align}\label{3d_DNLSeq:4.42}
\|{\bm \eta}^{m}\|_{h}^{2}+\tau|\widehat{\bm \eta}^{m-1}|_{1,h}^{2}\leq \widehat{C}_{0}(N^{-s}+\tau^{2})^{2},
\end{align}
for $2\leq m\leq M$. Clearly, it follows from Lemma \ref{3d_DNLSlem4.6} that \eqref{3d_DNLSeq:4.42} holds for $m=2$.
By the mathematical induction, we suppose that \eqref{3d_DNLSeq:4.42} holds at the first $n$-th step, that is, $m\leq n$, and we need to find such a $\widehat{C}_{0}$, independent of $n, \tau$ and $h$, that \eqref{3d_DNLSeq:4.42} holds  for $m=n+1$.

We multiply \eqref{3d_DNLS:4.16} by $2\widehat{\eta}_{\vec{j}}^{n}$ and sum them up for $\vec{j}\in J_{h}$ to arrive at
\begin{align}\label{3d_DNLSeq:4.34}
\text{i}\delta_{t}\|{\bm \eta}^{n}\|_{h}^{2}-\frac{1}{\tau}\text{Im}\langle {\bm \eta}^{n+1},{\bm \eta}^{n-1}\rangle_{h}- 2|\widehat{\bm \eta}^{n}|_{h}^{2}+R^{n}=\langle{\bm \xi}^{n},2\widehat{\bm \eta}^{n}\rangle_{h},
\end{align}
where
\begin{align*}
R^{n}=2h_{\Delta}\sum_{\vec{j}\in J_{h}}\Big[\widehat{({f}^{*})}_{\vec{j}}^{n}-\beta e^{-2\gamma t_{n}}|{U}_{\vec{j}}^{n}|^{2}
\widehat{U}_{\vec{j}}^{n}\Big]\widehat{\bar{\eta}}_{\vec{j}}^{n}.
\end{align*}
$R^{n}$ can be written as
\begin{align*}
R^{n}&=2\langle \widehat{{\bm f}^{*}}^{n}-\widehat{\bm f}^{n}, \widehat{\bm \eta}^{n}\rangle_{h}+ 2\langle \widehat{\bm f}^{n}-F({\bm u}^{n}), \widehat{\bm \eta}^{n}\rangle_{h}\\
&+2\langle F({\bm u}^{n})-F(({\bm u}^{*})^{n}), \widehat{\bm \eta}^{n}\rangle_{h}+
2\langle F(({\bm u}^{*})^{n})-F({\bm U}^{n}), \widehat{\bm \eta}^{n}\rangle_{h}\nonumber\\
&:=R_{1}^{n}+R_{2}^{n}+R_{3}^{n}+R_{4}^{n},
\end{align*}
where
\begin{align*}
F({\bm v}^n)=\beta e^{-2\gamma t_{n}}|{\bm v}^n|^{2}\widehat{\bm v}^n.
\end{align*}
From the Taylor formula and Lemma \ref{3d_DNLSlem4.3}, we can derive
\begin{align}\label{3d_DNLSeq:4.44}
R_{1}^{n}\leq CN^{-2s}+\|\widehat{\bm \eta}^{n}\|_{h}^{2},\ R_{2}^{n}\leq C\tau^{4}+||\widehat{\bm \eta}^{n}||_{h}^{2}.
\end{align}
According to Lemma \ref{3d_DNLSlem4.2}, it holds that
\begin{align}\label{3d_DNLSeq:4.45}
\|({\bm u}^{*})^{n}\|_{h,\infty}&\leq \max\limits_{{\bm x}\in \Omega}|u^{*}(t_{n},{\bm x})|\leq C\|u^{*}(t_{n},{\bm x})\|_{2}\leq C\|u(t_{n},{\bm x})\|_{2}\leq C,
\end{align}
 for any $n\ge 0$.
Then, with Lemma \ref{3d_DNLSlem4.3}, we can obtain
\begin{align}\label{3d_DNLSeq:4.46}
R_{3}^{n}\leq CN^{-2s}+\|\widehat{\bm \eta}^{n}\|_{h}^{2}.
\end{align}
We rewrite $R_{4}^{n}$ by
\begin{align}\label{3d_DNLSeq:4.47}
R_{4}^{n}&=2h_{\Delta}\sum_{\vec{j}\in J_{h}}\Big[F((u^{*})_{\vec{j}}^{n})-F(U_{\vec{j}}^{n})\Big]\widehat{\bar{\eta}}_{\vec{j}}^{n}\nonumber\\
&=2\beta e^{-2\gamma t_{n}} h_{\Delta}\sum_{\vec{j}\in J_{h}}\Big[|(u^{*})_{\vec{j}}^{n}|^{2}\widehat{(u^{*})}_{\vec{j}}^{n}
-|U_{\vec{j}}^{n}|^{2}\widehat{U}_{\vec{j}}^{n}\Big]\widehat{\bar{\eta}}_{\vec{j}}^{n}\nonumber\\
&=2\beta e^{-2\gamma t_{n}} h_{\Delta}\sum_{\vec{j}\in J_{h}}|(u^{*})_{\vec{j}}^{n}|^{2}|\widehat{\eta}_{\vec{j}}^{n}|^{2}\nonumber\\
&~~~+2\beta e^{-2\gamma t_{n}} h_{\Delta}\sum_{\vec{j}\in J_{h}}\Big(|(u^{*})_{\vec{j}}^{n}|^{2}-|U_{\vec{j}}^{n}|^{2}\Big)\Big(\widehat{(u^{*})}_{\vec{j}}^{n}
-\widehat{\eta}_{\vec{j}}^{n}\Big)\widehat{\bar{\eta}}_{\vec{j}}^{n}\nonumber\\
&:=R_{41}^{n}+R_{42}^{n}.
\end{align}
We then can deduce that
\begin{align}\label{3d_DNLSeq:4.48}
\text{Im}\big(R_{41}^{n}\big)=0,
\end{align}
and
\begin{align}\label{3d_DNLSeq:4.49}
|R_{42}^{n}|&\leq 2|\beta| e^{-2\gamma t_{n}} h_{\Delta}\sum_{\vec{j}\in J_{h}}\Big[|\eta_{\vec{j}}^{n}|^2+2|(u^{*})_{\vec{j}}^{n}|\cdot |\eta_{\vec{j}}^{n}|\Big]\Big(\Big|\widehat{(u^{*})}_{\vec{j}}^{n}\Big|
+\Big|\widehat{\eta}_{\vec{j}}^{n}\Big|\Big)\Big|\widehat{\eta}_{\vec{j}}^{n}\Big|\nonumber\\
&\le2|\beta| e^{-2\gamma t_{n}} h_{\Delta}\sum_{\vec{j}\in J_{h}}\Big[|\eta_{\vec{j}}^{n}|^2\Big|\widehat{\eta}_{\vec{j}}^{n}\Big|^2+
|\eta_{\vec{j}}^{n}|^2\Big|\widehat{(u^*)}_{\vec{j}}^{n}\Big|\Big|\widehat{\eta}_{\vec{j}}^{n}\Big|\nonumber\\
&+2\Big|(u^*)_{\vec{j}}^{n}\Big||\eta_{\vec{j}}^{n}|\Big|\widehat{\eta}_{\vec{j}}^{n}\Big|^2
+2\Big|(u^*)_{\vec{j}}^{n}\Big||\eta_{\vec{j}}^{n}|\Big|\widehat{(u^*)}_{\vec{j}}^{n}\Big|\Big|\widehat{\eta}_{\vec{j}}^{n}\Big|\nonumber\\
&\le Ch_{\Delta}\sum_{\vec{j}\in J_{h}}\Big[|\eta_{\vec{j}}^{n}|^2\Big|\widehat{\eta}_{\vec{j}}^{n}\Big|^2+|\eta_{\vec{j}}^n|^2+\Big|\widehat{\eta}_{\vec{j}}^{n}\Big|^2\Big]\nonumber\\
&\leq C\|{\bm \eta}^{n}\|_{h}^{2}+C\|\widehat{\bm \eta}^{n}\|_{h}^{2}+C\|{\bm \eta}^{n}\|_{h,3}^{2}\|\widehat{\bm \eta}^{n}\|_{h,6}^{2},
\end{align}
where the H\"{o}lder inequality is used.
With the above inequalities, we can prove that
\begin{align}\label{3d_DNLSeq:4.50}
|\text{Im}(R^n)|\leq C\|{\bm \eta}^{n}\|_{h}^{2}+C\|\widehat{\bm \eta}^{n}\|_{h}^{2}+C\|{\bm \eta}^{n}\|_{h,3}^{2}\|\widehat{\bm \eta}^{n}\|_{h,6}^{2}+C(N^{-s}+\tau^{2})^{2}.
\end{align}
The imaginary part of \eqref{3d_DNLSeq:4.34} yields
\begin{align}\label{3d_DNLSeq:4.51}
\frac{1}{2\tau}\Big(\|{\bm \eta}^{n+1}\|_{h}^{2}-\|{\bm \eta}^{n-1}\|_{h}^{2}\Big)+\text{Im}(R^{n})=\text{Im}\langle{\bm \xi}^{n},2\widehat{\bm \eta}^{n}\rangle_{h}.
\end{align}
According to the Cauchy-Schwartz inequality and \eqref{3d_DNLSeq:4.50}, \eqref{3d_DNLSeq:4.51} reduces to
\begin{align}\label{3d_DNLSeq:4.52}
\frac{1}{\tau}\Big(\|{\bm \eta}^{n+1}\|_{h}^{2}-\|{\bm \eta}^{n-1}\|_{h}^{2}\Big)&\le C\|{\bm \eta}^{n}\|_{h}^{2}+C\|\widehat{\bm \eta}^{n}\|_{h}^{2}\nonumber\\
&+C\|{\bm \eta}^{n}\|_{h,3}^{2}\|\widehat{\bm \eta}^{n}\|_{h,6}^{2}+C(N^{-s}+\tau^{2})^{2}.
\end{align}
The real part of \eqref{3d_DNLSeq:4.34} yields
\begin{align}\label{3d_DNLSeq:4.53}
|\widehat{\bm \eta}^{n}|_{h}^{2}+\frac{1}{2\tau}\text{Im}\langle {\bm \eta}^{n+1},{\bm \eta}^{n}\rangle_{h}-\frac{1}{2}\text{Re}\big(R^{n}\big)=-\text{Re}\langle{\bm \xi}^{n},\widehat{\bm \eta}^{n}\rangle_{h}.
\end{align}
Thanks to \eqref{3d_DNLSeq:4.45}, it follows from \eqref{3d_DNLSeq:4.47} that
\begin{align}\label{3d_DNLSeq:4.54}
|R_{41}^{n}|\leq C\|\widehat{\bm \eta}^{n}\|^{2}.
\end{align}
With Lemma \ref{3d_DNLSlem2.2}, the Cauchy-Schwartz inequality, \eqref{3d_DNLSeq:4.49} and \eqref{3d_DNLSeq:4.54}, \eqref{3d_DNLSeq:4.53} reduces to
\begin{align}\label{3d_DNLSeq:4.55}
|\widehat{\bm \eta}^{n}|_{1,h}^{2}&\leq C\|{\bm \eta}^{n}\|_{h}^{2}+C\|\widehat{\bm \eta}^{n}\|_{h}^{2}+C\|{\bm \eta}^{n}\|_{h,3}^{2}\|\widehat{\bm \eta}^{n}\|_{h,6}^{2}\nonumber\\
&+\frac{1}{4\tau}\Big(\|{\bm \eta}^{n+1}\|_{h}^{2}+\|{\bm \eta}^{n-1}\|_{h}^{2}\Big)+C(N^{-s}+\tau^{2})^{2}.
\end{align}

To reduce the nonlinear term in the above inequality, now we prove the following inequality
\begin{align}\label{3d_DNLSeq:4.56}
||{\bm \eta}^{n}||_{h,3}^{2}||\widehat{\bm \eta}^{n}||_{h,6}^{2}&\leq C\big( ||{\bm \eta}^{n+1}||_{h}^{2}+||{\bm \eta}^{n}||_{h}^{2}+||{\bm \eta}^{n-1}||_{h}^{2} \big)\nonumber\\
&+C||\widehat{\bm \eta}^{n}||_{h}^{2}+C(N^{-s}+\tau^{2})^{2},
\end{align}
with two different cases.

Firstly, we consider the case $\tau\leq h\ (h=\frac{2\pi}{N})$. We use Lemma \ref{3d_DNLSlem2.3} and \eqref{3d_DNLSeq:4.42} to get
\begin{align*}
&\|{\bm\eta}^{n}\|_{h,3}\leq Ch^{-\frac{1}{2}}\|{\bm \eta}^{n}\|_{h}\leq Ch^{\frac{3}{2}},\ \|\widehat{\bm\eta}^{n}\|_{6,h}\leq Ch^{-1}\|\widehat{\bm \eta}^{n}\|_{h}.
\end{align*}
When $ Ch\leq 1,$
\begin{align}\label{3d_DNLSeq:4.59}
\|{\bm \eta}^{n}\|_{h,3}^{2}\|\widehat{\bm \eta}^{n}\|_{h,6}^{2}\leq C\|\widehat{\bm \eta}^{n}\|_{h}^{2}.
\end{align}

Secondly, we consider the case $\tau\ge h$. By Lemma \ref{3d_DNLSlem2.3}, we get
\begin{align}\label{3d_DNLSeq:4.60}
\|\widehat{\bm \eta}^{m}\|_{h,6}^{2}\leq C\big(|\widehat{\bm \eta}^{m}|_{1,h}+\|\widehat{\bm \eta}^{m}\|_{h})^{2}.
\end{align}
By using \eqref{3d_DNLSeq:4.42}, together with Lemma \ref{3d_DNLSlem2.4} and the H\"older inequality, we have
\begin{align}\label{3d_DNLSeq:4.61}
\|{\bm \eta}^{m}\|_{h,3}^{2}&\leq \|{\bm \eta}^{m}\|_{h}\|{\bm \eta}^{m}\|_{h,6}\nonumber\\
&\leq \|{\bm \eta}^{m}\|_{h}\big(|{\bm \eta}^{m}|_{1,h}+\|{\bm \eta}^{m}\|_{h})\nonumber\\
&\leq \|{\bm \eta}^{m}\|_{h}\Big(\frac{2T}{\tau}\max\limits_{1 \leq k\leq
m-1}|\widehat{\bm \eta}^{k}|_{1,h}+|{\bm \eta}^{1}|_{1,h}+|{\bm \eta}^{0}|_{1,h}+\|{\bm \eta}^{m}\|_{h}\Big)\nonumber\\
&\leq C\tau^{\frac{5}{2}},
\end{align}
for $1\leq m\leq n$.

Combining \eqref{3d_DNLSeq:4.60} with \eqref{3d_DNLSeq:4.61}, we obtain
\begin{align}\label{3d_DNLSeq:4.62}
\|{\bm \eta}^{n}\|_{h,3}^{2}\|\widehat{\bm \eta}^{n}\|_{h,6}^{2}&\leq \epsilon\tau\big(|\widehat{\bm \eta}^{n}|_{1,h}^{2}+\|\widehat{\bm \eta}^{n}\|_{h}^{2}\big),
\end{align}
when $C\tau^{\frac{3}{2}}\leq \epsilon$.
With the use of \eqref{3d_DNLSeq:4.62}, we can deduce from \eqref{3d_DNLSeq:4.55} that
\begin{align*}
|\widehat{\bm \eta}^{n}|_{1,h}^{2}&\leq C\|{\bm \eta}^{n}\|_{h}^{2}+C\|\widehat{\bm \eta}^{n}\|_{h}^{2}+\epsilon\tau|\widehat{\bm \eta}^{n}|_{1,h}^{2}\nonumber\\
&+\frac{1}{4\tau}\Big(\|{\bm \eta}^{n+1}\|_{h}^{2}+\|{\bm \eta}^{n-1}\|_{h}^{2}\Big)+C(N^{-s}+\tau^{2})^{2},
\end{align*}
which further implies that
\begin{align}\label{3d_DNLSeq:4.63}
|\widehat{\bm \eta}^{n}|_{1,h}^{2}&\leq C\|{\bm \eta}^{n}\|_{h}^{2}+C\|\widehat{\bm \eta}^{n}\|_{h}^{2}\nonumber\\
&+\frac{1}{4\tau}\Big(\|{\bm \eta}^{n+1}\|_{h}^{2}+\|{\bm \eta}^{n-1}\|_{h}^{2}\Big)+C(N^{-s}+\tau^{2})^{2},
\end{align}
when $\tau$ sufficiently small.
With \eqref{3d_DNLSeq:4.63}, \eqref{3d_DNLSeq:4.62} reduces to
\begin{align}\label{3d_DNLSeq:4.64}
\|{\bm \eta}^{n}\|_{h,3}^{2}\|\widehat{\bm \eta}^{n}\|_{h,6}^{2}&\leq C\big(\|{\bm \eta}^{n+1}\|_{h}^{2}+\|{\bm \eta}^{n}\|_{h}^{2}+\|{\bm \eta}^{n-1}\|_{h}^{2} \big)\nonumber\\
&+C\|\widehat{\bm \eta}^{n}\|_{h}^{2}+C(N^{-s}+\tau^{2})^{2}.
\end{align}
Up to now, we have proved that \eqref{3d_DNLSeq:4.56} holds.

By using \eqref{3d_DNLSeq:4.56}, we can deduce from \eqref{3d_DNLSeq:4.52} that
\begin{align}\label{3d_DNLSeq:4.65}
\frac{1}{\tau}\Big(\|{\bm \eta}^{n+1}\|_{h}^{2}-\|{\bm \eta}^{n-1}\|_{h}^{2}\Big)&\le C\big(\|{\bm \eta}^{n+1}\|_{h}^{2}+\|{\bm \eta}^{n}\|_{h}^{2}+\|{\bm \eta}^{n-1}\|_{h}^{2}\big)\nonumber\\
&+C(N^{-s}+\tau^{2})^{2}.
\end{align}
Applying Gronwall's inequality \cite{zhou90} to \eqref{3d_DNLSeq:4.65}, we obtain
\begin{align}\label{3d_DNLSeq:4.66}
\|{\bm \eta}^{n+1}\|_{h}^{2}\le CTe^{2CT}(N^{-s}+\tau^{2})^{2},
\end{align}
where $\tau$ is sufficiently small, such that $C\tau\le \frac{1}{2}$.
Moreover, from \eqref{3d_DNLSeq:4.55}-\eqref{3d_DNLSeq:4.56} and \eqref{3d_DNLSeq:4.66}, we see
\begin{align}\label{3d_DNLSeq:4.67}
\tau|\widehat{\bm \eta}^{n}|_{1,h}^{2}\leq CTe^{2CT}(N^{-s}+\tau^{2})^{2}.
\end{align}
Thus, \eqref{3d_DNLSeq:4.42} holds for $m=n+1$, if we take $\widehat{C}_{0}=2CTe^{2CT}$. We complete the induction.

According to Lemma \ref{3d_DNLSlem4.3} and \eqref{3d_DNLSeq:4.66}, we can deduce
\begin{align}\label{eq:4.68}
\|{\bm u}^{n}-{\bm U}^{n}\|_{h}\leq \|{\bm u}^{n}-({\bm u}^{*})^{n}\|_{h}+\|{\bm \eta}^{n}\|_{h}\leq C(N^{-s}+\tau^{2}).
\end{align}
which further implies that
\begin{align*}
\|{\bm \psi}^{n}-{\bm \Psi}^{n}\|_{h}\leq  Ce^{-\gamma t_{n}}(N^{-s}+\tau^{2}).
\end{align*}
This completes the proof.

\begin{rmk}\rm
When $\gamma=0$, the proposed scheme reduces to a linearly implicit and conservative Fourier pseudo-spectral scheme of the NLS equation in 3D and the error estimate is also valid. In addition,
the analysis techniques used in this paper can be extended to establish an optimal $L^2$-error
estimate for the NLS equation with the cubic damped term \cite{BJM04}.
\end{rmk}

\section{Numerical experiments}
In this section, we will investigate the numerical behaviors of the scheme \eqref{3d_NLS:eq2.3}-\eqref{3d_NLS:eq2.4}.
Also, the results are compared
with the IFD scheme \cite{zhang04} the classical Runge-Kutta method of order 3 (denoted by RK3 method) for numerical errors and CPU times, respectively.
{For clarity, two selected schemes are given, respectively, as follows:
\begin{itemize}
\item IFD scheme:  \begin{align}\label{IFD-scheme}
&\text{i} \delta_t^+{\bm U}^{n}+{\bm\Delta}_{1,h}{\bm U}^{n+\frac{1}{2}}+\frac{\beta}{2}e^{-2\gamma t_{n+\frac{1}{2}}}\Big(|{\bm U}^{n+1}|^{2}
+|{\bm U}^{n}|^{2}\Big)\cdot{\bm U}^{n+\frac{1}{2}}=0,
\end{align}
where ${\bm\Delta}_{1,h}={\bm I}_{N_{3}}\otimes {\bm I}_{N_{2}}\otimes{\bm B}_{1}+{\bm I}_{N_{3}}\otimes {\bm B}_{2}\otimes {\bm I}_{N_{1}}
+{\bm B}_{3}\otimes {\bm I}_{N_{2}}\otimes {\bm I}_{N_{1}}$ represents the conventional finite difference method discretization (see \eqref{finite-difference-method}) of the Laplace operator.
\item RK3 method:
\begin{align}\label{3rd-RK}
\left \{
 \aligned
&{\bm K}_1=f(t_n,{\bm U}^n),\\
&{\bm K}_2=f(t_n+\frac{\tau}{2},{\bm U}^n+\frac{\tau}{2}{\bm K}_1),\\
&{\bm K}_3=f(t_n+\tau,{\bm U}^n-\tau{\bm K}_1+2\tau{\bm K}_2),\\
&{\bm U}^{n+1}={\bm U}^n+\frac{\tau}{6}({\bm K}_1+4{\bm K}_2+{\bm K}_3),
\endaligned
 \right.
\end{align}
where $f(t,{\bm U})=\text{i}\big({\bm\Delta}_{h}{\bm U}+\beta e^{-2\gamma t}|{\bm U}|^2\cdot {\bm U}\big)$ and ${\bm\Delta}_{h}$ represents the Fourier pseudo-spectral method discretization (see \eqref{Laplace-discrete}) of the Laplace operator.
\end{itemize}}
\noindent As a summary, a detailed table on the properties of each scheme has been
given in Table. \ref{Tab_2DNLS:1}. {Note that the IFD scheme \eqref{IFD-scheme} satisfies the discrete mass conservation law \eqref{mass-conservation-law} and the following discrete energy conservation law
\begin{align}\label{IFD-energy-conservation-law}
\mathcal{E}^{n}=\mathcal{E}^{0},\ n=1,\cdots,M,
\end{align}
where
\begin{align*}
\mathcal{E}^{n}=&|{\bm U}^{n}|_{1,h}^{2}-\frac{\beta}{2}e^{-2\gamma t_{n-\frac{1}{2}}}\|{\bm U}^{n}\|_{h,4}^{4}-\frac{\beta}{2}\sum_{l=1}^{n}e^{-2\gamma t_{l-\frac{3}{2}}}(1-e^{-2\gamma\tau})\|{\bm U}^{l-1}\|_{h,4}^{4}.
\end{align*}}

{ In order to quantify the residuals of the discrete conservation laws, we use the relative mass residual and the relative energy residual between the discrete mass (i.e., $\mathcal{M}^n$) as well as the discrete energy (i.e., $\mathcal{E}^n$) at $t=t_n$ and the initial discrete ones, respectively, as
\begin{align}
RM^n=\Big|\frac{\mathcal{M}^n-\mathcal{M}^0}{\mathcal{M}^0}\Big|,\ RE^n=\Big|\frac{\mathcal{E}^n-\mathcal{E}^0}{\mathcal{E}^0}\Big|,\ 0\le n\le M.
\end{align}}

In our computations, the Jacobi iteration method is used to solve the LI-CFP scheme \eqref{3d_NLS:eq2.3} with the tolerance number $10^{-14}$,
while, for the IFD scheme \eqref{IFD-scheme}, we use the following fixed-point iteration method to solve all the nonlinear algebraic equations
\begin{align*}
&\text{i} \frac{{\bm U}^{n+\frac{1}{2},(s+1)}-{\bm U}^{n}}{\tau}+\frac{1}{2} {\bm\Delta}_{1,h}{\bm U}^{n+\frac{1}{2},(s+1)}\nonumber\\
&~~~~~~~~~~+\frac{\beta}{4}e^{-2\gamma t_{n+\frac{1}{2}}}\Big[\Big(|2{\bm U}^{n+\frac{1}{2},(s)}-{\bm U}^{n}|^{2}
+|{\bm U}^{n}|^{2}\Big)\cdot{\bm U}^{n+\frac{1}{2},(s)}\Big]=0,
\end{align*}
where ${\bm  U}^{n+\frac{1}{2}}=\frac{{\bm  U}^{n+1}+{\bm  U}^{n}}{2}$. We set $10^{-14}$ as the error tolerance and the linear system is solved efficiently by the fast solver presented in Ref. \cite{JCWL17} in every iteration. Solving the above equations gives ${\bm U}^{n+\frac{1}{2}}$. Then, we have ${\bm U}^{n+1}=2{\bm  U}^{n+\frac{1}{2}}-{\bm U}^{n}$.

To evaluate the convergence rate, for a fixed $n$, we use the formula
\begin{align*}
\text{Rate}=\frac{\ln( error_{1}/error_{2})}{\ln (\tau_{1}/\tau_{2})},
\end{align*}
where $\tau_{\tilde{l}}, error_{\tilde{l}}, (\tilde{l}=1,2)$ are step sizes and errors with the step size $\tau_{\tilde{l}}$, respectively.

\begin{table}[H]
\tabcolsep=9pt
\small
\renewcommand\arraystretch{1.1}
\centering
\caption{Comparison of properties of different numerical schemes}\label{Tab_2DNLS:1}
\begin{tabular*}{\textwidth}[h]{@{\extracolsep{\fill}}c c c c c c}\hline 
 \diagbox{Property}{Scheme/Method}& LI-CFP scheme& IFD scheme& RK3 method\\\hline
 Mass conservation& Yes&Yes&No \\[1ex]
 Energy conservation& Yes&Yes&No \\[1ex]
 Fully implicit&No&Yes&No\\[1ex]
 Linearly implicit&Yes&No&No\\[1ex]
 Fully explicit&No&No&Yes\\[1ex]
 Temporal accuracy& 2nd&2nd&3rd \\[1ex]
\hline
\end{tabular*}
\end{table}
{ Eq. \eqref{3DNLD_eq:1.1} possesses the following analytical solution
\begin{align}\label{3d_NLS:eq6.1}
\psi(t,{\bm x})=Ke^{-\gamma t}e^{\text{i}(k_1 x+k_2y+k_3z-\delta(t))},\ {\bm x}\in\Omega=[0,2\pi]^3,
\end{align}
where $\delta(t)=\big(k_1^2+k_2^2+k_3^2\big)t+\frac{\beta|K|^2}{2\gamma}(e^{-2\gamma t}-1)$.}

We choose $K=k_1=k_2=k_3=1,\beta=2$ and set the analytical solution of \eqref{3d_NLS:eq6.1} at $t=0$ on the domain $\Omega=[0,2\pi]^3$ as the
initial condition. The convergence rates at $t=1$ in the temporal direction of the proposed scheme are displayed in Table. \ref{3d-DNLS-Tab1}. As illustrated in the table, the convergence order of our scheme is of second order accuracy for $\gamma=1$ in temporal
direction, which verifies the error estimate in Theorem \ref{3d_DNLSthm4.1}. We then depict the spatial errors of the LI-CFP scheme with $\gamma=1$ and $\tau=10^{-5}$ at $t=1$ in Table. \ref{3d-DNLS-Tab2}, which shows that the spatial error of the proposed scheme is very small and almost negligible, and the error is dominated
by the time discretization error. It confirms that, for sufficiently smooth problems,
the Fourier pseudo-spectral method is of arbitrary order in space. Finally, we list the errors and CPU times of the three different schemes (see Table \ref{Tab_2DNLS:1}) in the solutions at $t=1$ in Table. \ref{Tab521}.  From the table, we can draw two main observations: (i) the RK3 method \eqref{3rd-RK} provides more accuracy and more efficient than the two schemes for short time computation; (ii) compared with the IFD scheme \eqref{IFD-scheme}, our scheme admits smaller numerical errors and more efficient for a fixed temporal and spatial step. {Here, we should note that the Tables \ref{3d-DNLS-Tab1}, \ref{3d-DNLS-Tab2}, and \ref{Tab521} show the error between the numerical solution ${\bm \Psi}^n\in\mathbb{V}_h$
and the exact solution ${\bm \psi}^n\in\mathbb{V}_h$ for a given time $t$.}

\begin{figure}[H]
\centering\begin{minipage}[t]{60mm}
\includegraphics[width=65mm]{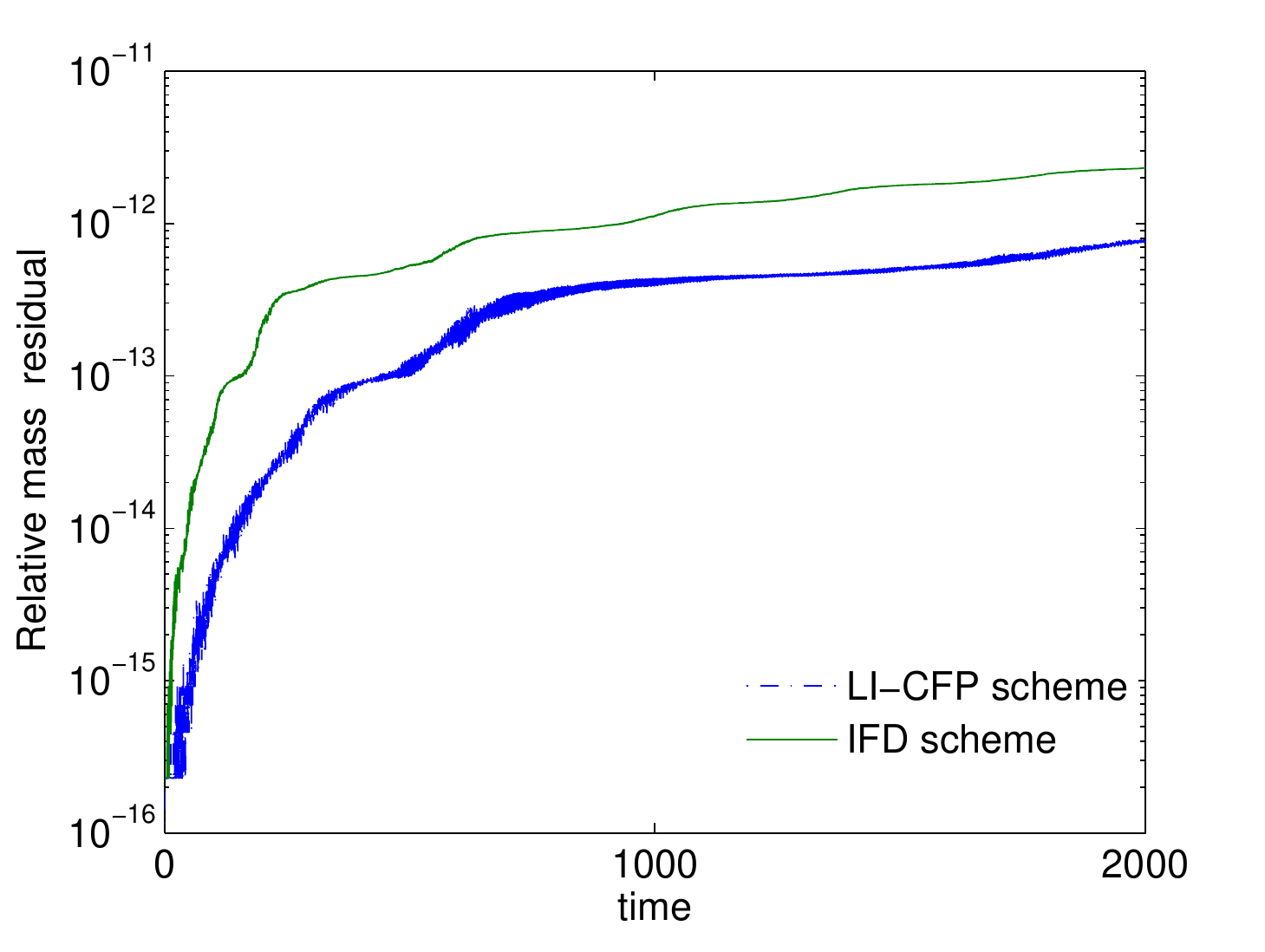}
\end{minipage}\ \
\begin{minipage}[t]{60mm}
\includegraphics[width=65mm]{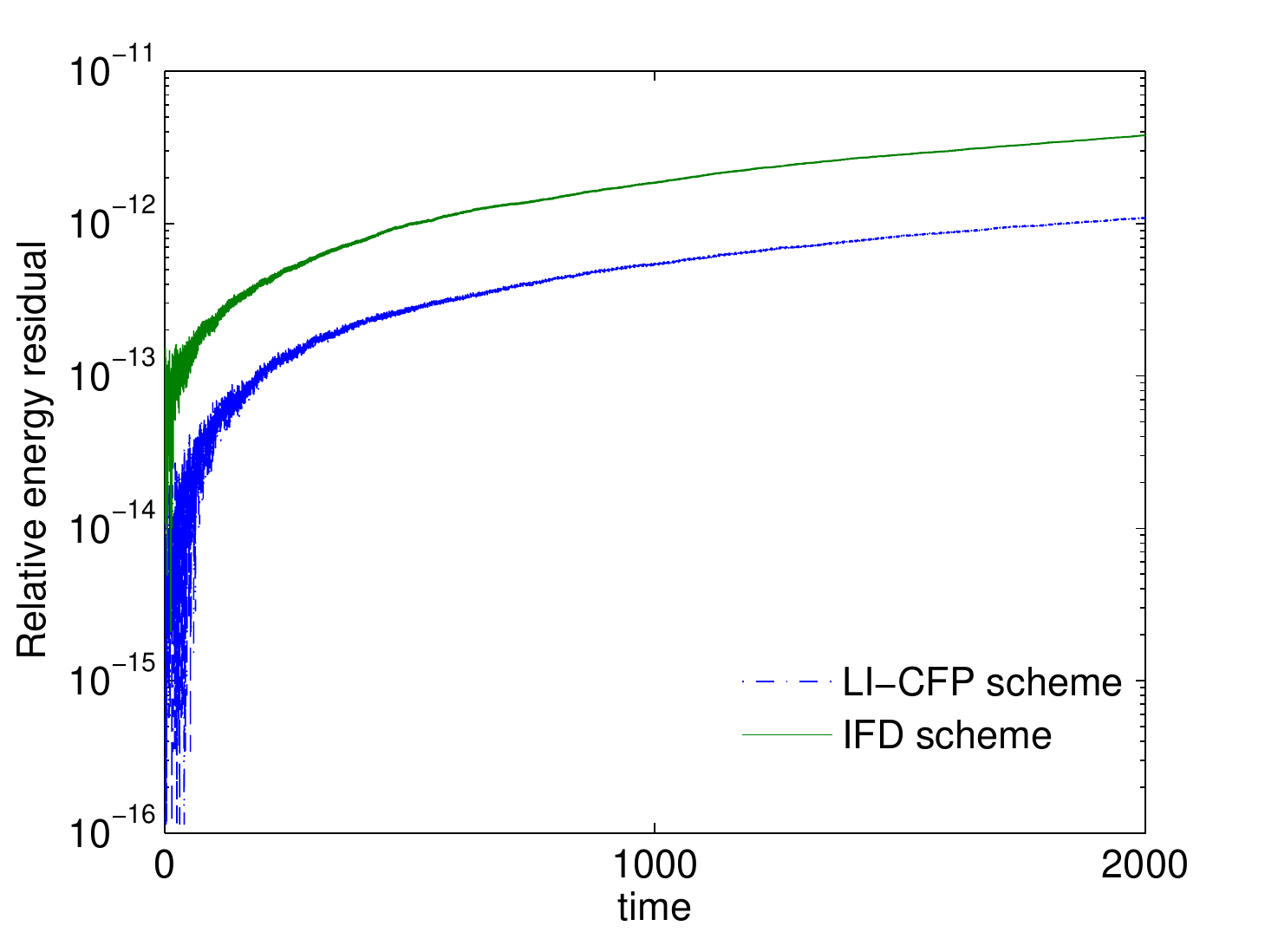}\\
\end{minipage}\\
\caption{The {relative residuals} in the mass (left) and the energy (right) over the time interval $t\in[0,2000]$ with $\tau=0.1$, $\gamma=1$ and $N_1=N_2=N_3=16$.}\label{fig523}
\end{figure}

To further investigate the discrete conservation laws of the proposed scheme and the IFD scheme, we provide the relative mass and energy residuals of the two schemes over the time interval $t\in[0,2000]$ in Figure. \ref{fig523}, which shows that both schemes can preserve the mass and energy conservation laws exactly,
 and the LI-CFP scheme admits smaller residuals than the one provided by the IFD scheme. 

\begin{table}[H]
\tabcolsep=9pt
\small
\renewcommand\arraystretch{1.2}
\centering
\caption{The temporal error and convergence rate of the proposed scheme with $\gamma=1$ and $N_1=N_2=N_3=16$ at $t=1$.}\label{3d-DNLS-Tab1}
\begin{tabularx}{\textwidth}{XXXXXX}\hline
{$\tau$}  &{$L^{2}$-norm} & {Rate}  & {$L^{\infty}$-norm} & {Rate}\\     
\hline
{0.1}  & {1.751e-01} & {-} & {1.112e-02}&{-}\\[1ex]
 {0.05} &  {4.507e-02}  & {1.958} & {2.862e-03}&{1.958}  \\[1ex]   
{0.025}  & {1.135e-02} & {1.989} & {7.209e-04}&{1.989} \\[1ex]
{0.0125} & {2.844e-03} & {1.997}  & {1.806e-04} & {1.997}  \\[1ex]
\hline
\end{tabularx}
\end{table}

\begin{table}[H]
\tabcolsep=9pt
\small
\renewcommand\arraystretch{1.2}
\centering
\caption{The spatial error of the proposed scheme with $\gamma=1$ and $\tau=10^{-5}$ at $t=1$.}\label{3d-DNLS-Tab2}
\begin{tabularx}{1.0\textwidth}{XXXX}\hline
{$N_1\times N_2\times N_3$}  &{$L^{2}$-norm}  & {$L^{\infty}$-norm} \\     
\hline
 {$4\times 4\times 4$}  & {1.796e-09}  & {1.141e-010}\\[1ex]
{$8\times 8\times 8$} &  {1.868e-09}   & {1.188e-010}  \\[1ex]   
{$16\times 16\times 16$}  & {1.818e-09}  & {1.162e-010} \\[1ex]
\hline
\end{tabularx}
\end{table}

\begin{table}[H]
\tabcolsep=13pt
\small
\renewcommand\arraystretch{2}
\centering
{\caption{The numerical errors and CPU times for different schemes with $\gamma=1$ at $t=1$}\label{Tab521}
\setlength{\tabcolsep}{2mm}{
\begin{tabularx}{\textwidth}{XXXXXXXX}\hline
{$\tau/h$ }&{Scheme/Method} &{$L^{2}$-norm}& {$L^{\infty}$-norm}&{CPU (s)}\\
\hline
{$0.005/\frac{\pi}{4}$}& {LI-CFP scheme}  &{4.553e-04}&{2.891e-05}&{0.9}\\[1ex]
{}   & {IFD scheme}  &{8.746e-01}&{5.553e-02}&{1.7}\\[1ex]
{}   & {RK3 method}  &{1.014e-06}&{6.441e-08}&{0.6} \\[1ex]
{$0.0025/\frac{\pi}{8}$}& {LI-CFP scheme}   &{1.138e-04}&{7.227e-06}&{8.4}\\[1ex]
{}   & {IFD scheme}  &{2.222e-01}&{1.411e-02}&{18.3}\\[1ex]
{}   & {RK3 method}  &{1.266e-07}&{8.041e-09}&{4.2} \\[1ex]
{$0.00125/\frac{\pi}{16}$}& {LI-CFP scheme}&{2.846e-05}&{1.807e-06}&{116.7}\\[1ex]
{}   & {IFD scheme}  &{5.578e-02}&{3.542e-03}&{272.3}\\[1ex]
{}   & {RK3 method}  &{1.582e-08}&{1.004e-09}&{75.4} \\[1ex]
\hline
\end{tabularx}}}
\end{table}

\section{Concluding remarks}
In this paper, we propose a linearly implicit and conservative Fourier pseudo-spectral scheme for the DNLS equation in 3D. We show that the proposed scheme is uniquely solvable, and preserves both mass and energy conservation laws.
We first introduce the semi-norm equivalence between the Fourier pseudo-spectral method and the finite difference
method, in order to establish an optimal error estimate of the proposed scheme. Then, we prove that, without
any restriction on the grid ratio, the proposed pseudo-spectral scheme is convergent with order $O(N^{-s}+\tau^2)$  in discrete $L^2$-norm.
Finally, numerical results verify the theoretical analysis.

{We conclude this paper with some remarks. First, compared with the IFD scheme, our scheme is more efficient and has a significant advantage in preserving the discrete conservation laws. Second, the RK3 method is more accurate and more efficient than the proposed method, however, such method is conditionally stable. Finally, to the best of our knowledge, the construction and numerical analysis of higher order linearly implicit structure-preserving schemes are still not available for the DNLS equation \eqref{3DNLD_eq:1.1}, which is an interesting topic for future studies.}

\section*{Acknowledgments}
The authors would like to express sincere gratitude to the referees for their insightful
comments and suggestions. This work is supported by the National Natural Science Foundation of China (Grant Nos. 11771213, 11901513),
the National Key Research and Development Project of China (Grant Nos. 2018YFC0603500, 2018YFC1504205), the Yunnan Provincial Department of Education Science Research Fund Project (Grant No. 2019J0956) and the Science and Technology
Innovation Team on Applied Mathematics in Universities of Yunnan.

%

\end{document}